\documentclass[10pt]{article}

\usepackage{bm}
\usepackage{amsmath}
\usepackage{amsthm}
\usepackage{amsfonts}
\usepackage{amssymb}
\usepackage{amscd}
\usepackage{subeqnarray}
\usepackage{hyperref}
\usepackage{pstricks}
\usepackage{latexsym}
\usepackage{enumerate}
\usepackage{epsfig}
\usepackage{subfig}
\usepackage{graphicx}
\usepackage{framed,fancybox}
\usepackage{bbm}
\usepackage{multirow}
\usepackage[T1]{fontenc}
\usepackage{threeparttable}
\usepackage{rotating}
\usepackage{color,soul,xcolor}
\usepackage{physics}
\usepackage[version=4]{mhchem}
\usepackage{siunitx}
\usepackage{footnpag}			      	
\usepackage{longtable,tabularx}
\newcommand{\m}[1]{{\bf{#1}}}
\newcommand{\g}[1]{\boldsymbol #1}

\newcommand{\C}[1]{{\cal {#1}}}
\newcommand{\mbb}[1]{\mathbb{#1}}
\newcommand{\R}{\mbb{R}}
\newcommand{\T}{^{\sf T}}    

\newcommand{\lrsup}[3]{{\vphantom{#2}}^{#1}{{#2}}^{#3}}
\newcommand{\lsupddt}[2]{{\vphantom{\frac{d#2}{dt}}}^{#1}{\frac{d#2}{dt}}}

\makeatletter
\newcommand{\thickhline}{%
    \noalign {\ifnum 0=`}\fi \hrule height 1pt
    \futurelet \reserved@a \@xhline
}
\newcolumntype{"}{@{\hskip\tabcolsep\vrule width 1pt\hskip\tabcolsep}}
\makeatother

\long\def\symbolfootnote[#1]#2{\begingroup\def\thefootnote{\fnsymbol{footnote}}\footnote[#1]{#2}\endgroup}
\def\qed{\hfill{$\vcenter{\hrule height1pt \hbox{\vrule width1pt height5pt
    \kern5pt \vrule width1pt} \hrule height1pt}$} \medskip}
\title{\bf Minimum-Time Reorientation of Axisymmetric Rigid Spacecraft Using Three Controls }

\author{Elisha R.~Pager\thanks{Ph.D.~Candidate, Department of Mechanical and Aerospace Engineering, University of Florida, Gainesville, Florida 32611-6250. Email: epager@ufl.edu.} \\
Anil V.~Rao\thanks{Professor, Department of Mechanical and Aerospace Engineering, University of Florida, Gainesville, FL 32611-6250. E-mail:  anilvrao@ufl.edu.  Corresponding Author.} \vspace{12pt} \\ {\em University of Florida} \\ {\em Gainesville, FL 32611}
}

\date{}

\begin{document}
\maketitle
\thispagestyle{empty}

\begin{abstract}
A minimum-time reorientation of an axisymmetric rigid spacecraft controlled by three torques is studied. The orientation of the body is modeled such that the attitude kinematics are representative of a spin-stabilized spacecraft. The optimal control problem considered is shown to have a switching control structure.  Moreover, under certain assumptions, the solutions contain segments that lie on a singular arc. A numerical optimization study is performed using a recently developed method that is designed to accurately solve bang-bang and singular optimal control problems.  The optimality conditions for the resulting optimal control problem are derived and analyzed for a variety of cases. Also, the results obtained in this study are compared to a previous method existing in the literature. The key features of the optimized trajectories and controls are identified, and the aforementioned method for solving bang-bang and singular optimal control problems is shown to efficiently and accurately solve the problem under consideration.  
\end{abstract}

\renewcommand{\baselinestretch}{1}
\normalsize\normalfont

 \section{Nomenclature}

{\renewcommand\arraystretch{1.0}
   \noindent\begin{longtable*}{@{}l @{\quad=\quad} l@{}}
     $a$ & alternate inertial tensor constant \\
     $(\bf{b}_1,\bf{b}_2,\bf{b}_3)$ & principal-axis vector basis fixed in the body \\
     $\C{B}$ & rigid-body \\
     $C$ & center of mass of rigid-body $\C{B}$ \\
     $g$ & switching function \\
     $\C{H}$ & Hamiltonian \\
     $(I_1,I_2,I_3)$ & moments of inertia expressed in $\{C,\m{b} \}$\\
     $\m{I}^{\C{B}}$ & moment of inertia tensor of $\C{B}$ relative to $C$ \\
     $\C{J}$ & cost \\
     $\bf{M}$ & resulting moment \\
     $(\bf{n}_1,\bf{n}_2,\bf{n}_3)$ & vector basis fixed in the inertial reference frame \\
     $\C{N}$ & inertial reference frame \\
     $O$ & point fixed in inertial reference frame \\
     $p$ & regularization iteration index \\
     $r$ & index for number of differentiation's of switching function \\
     $t$ & time \\
     $t_f$ & terminal time \\
     $t_s$ & switch time \\
     $(u_1,u_2,u_3)$ & components of alternate control parameterization \\
     $u_{\min}$ & lower bound on control \\
     $u_{\max}$ & upper bound on control \\
     $u_s$ & singular control \\
     $(x_1,x_2)$ & relative position of the $\m{n}_3$ axis in the basis $\m{b}$ \\
     $\m{Y}$ & state vector \\
     $\g{\alpha}$ & angular acceleration vector \\
     $\delta$ & value of regularization term \\
     $\epsilon$ & regularization parameter \\
     $\epsilon_{\textrm{mesh}}$ & mesh refinement accuracy tolerance \\
     $\epsilon_{\textrm{NLP}}$ & NLP solver tolerance \\
     $\g{\lambda}$ & costate vector \\
     $\g{\nu}$ & constraint multiplier vector \\
     $(\g{\omega}_1,\g{\omega}_2,\g{\omega}_3)$ & components of angular velocity \\
     $\Phi$ & differentiable vector function that defines the terminal state conditions \\
     $(\sigma,\beta,\gamma)$ & direction cosines of $\m{n}_3$ w.r.t. $\m{b}$ \\
     $\g{\tau}$ & torque \\
     $\{C,\m{b} \}$ & principal-axis body-fixed coordinate system \\
     $\{O,\m{n} \}$ & inertially fixed coordinate system \\
   \end{longtable*}}

\renewcommand{\baselinestretch}{2}
\normalsize\normalfont 

\section{Introduction}

Minimum-time reorientation of a rigid body is a problem of great interest that has been studied for decades.  Such problems play an important role in spacecraft attitude control.  While different variations to minimum-time reorientation problems have been explored over the last few decades, numerical approaches and computation power have grown considerably allowing researchers to readdress these fundamental problems of interest using new tools and methods. In particular, the work of Ref.~\cite{Shen1999} performed an analysis for the optimal minimum-time reorientation of an axisymmetric rigid spacecraft under the influence of two control torques.  While most of the analysis of this fundamental problem remains relevant today, the numerical approaches used in these prior studies were developed many years ago.  In particular, solutions to the minimum-time reorientation problem may be such that the optimal control is either bang-bang or singular (that is, Pontryagin's minimum principle fails to yield a complete solution~\cite{Kirk2004,Bryson1975}), and these previously developed methods were not easily able to solve optimal control problems that contained such a switching structure.  

Minimum-time rigid-body reorientation problems first gained traction in the work of Ref.~\cite{Bilimoria1993} where a three-axis rigid body spacecraft was considered and eigenaxis rotations were shown to not be time optimal. The work of Ref.~\cite{Shen1999} later showed that more generic time optimal solutions could be obtained by considering axisymmetric rigid bodies controlled by two torques, i.e. an under actuated spacecraft. Together Refs.~\cite{Shen1999} and~\cite{Tsiotras1995} further advanced the work of Ref.~\cite{Bilimoria1993} by extending the analysis from a triaxisymmetric body to an axisymmetric body. This work highlighted the difficulties of solving these minimum-time reorientation problems using the methods available at that time. These difficulties included the need to produce accurate guesses for the costates and handling bang-bang and singular arcs in the control solution (specifically when using indirect methods). Despite these difficulties and limitations, the results of Refs.~\cite{Shen1999,Tsiotras1995} have practical applications, as many spacecraft can be approximated as an axisymmetric body instead of a triaxisymmetric rigid body. The results in Ref.~\cite{Fleming2008} later showed that the time-optimal reorientation of an axisymmetric body could be extended to three control torques while still considering a spin-stabilized orientation, and solutions obtained for this three-torque control problem showed improvement over solutions obtained obtained with the two-torque control model. More recently, this application has been further extended to consider three-axis rigid-bodies and the minimum-time reorientation of asymmetric rigid bodies as shown in Refs.~\cite{Kim2013,Li2016,Lan2020}.  In particular, Refs.~\cite{Kim2013,Li2016,Lan2020} assume singular arcs are not possible in the solutions and therefore such behavior is not considered in the analysis.

As seen from the prior literature, the complexity of minimum-time reorientation optimal control problems makes it necessary to solve such problems numerically.  Numerical methods for solving optimal control problems fall into two broad categories: indirect methods and direct methods. In an indirect method, the first-order variational optimality conditions are derived, and the optimal control problem is converted to a Hamiltonian boundary-value problem (HBVP). The HBVP is then solved numerically using a differential-algebraic equation solver. In a direct method, the control or state and control are approximated, and the optimal control problem is transcribed into a finite-dimensional nonlinear programming problem (NLP)~\cite{Betts2010}. The NLP is then solved numerically using well-developed software such as SNOPT~\cite{Gill2002} or IPOPT~\cite{Biegler2008}. Many of the earlier results to the minimum-time rigid-body reorientation problem were obtained using indirect methods while the more recent research has utilized direct methods.

Motivated by the ongoing interest in the optimal control of rigid bodies, this work revisits the minimum-time reorientation of a spin-stabilized axisymmetric rigid spacecraft under the influence of three control torques.  As stated earlier, the solutions to such problems are known to have an optimal control that is either bang-bang or singular \cite{Shen1999,Seywald1993}.  This work analyzes the optimal solution structure of the minimum-time reorientation problem and discusses optimal solutions of various minimum-time maneuvers. Furthermore, numerical solutions are obtained using a recently developed method, called the BBSOC method, for solving bang-bang and singular optimal control problems~\cite{Pager2022}.

This research is inspired by the work of Ref.~\cite{Shen1999}. Although both Ref.~\cite{Shen1999} and this work focus on the optimization of a minimum-time axisymmetric rigid spacecraft, the work in this paper differs significantly from the work of Ref.~\cite{Shen1999} in the following ways.  First, the numerical approach used to solve the problem in Ref.~\cite{Shen1999} consisted of using an indirect shooting method which required the formulation of the optimality conditions whereas in this work the aforementioned direct collocation BBSOC method~\cite{Pager2022} is employed. This difference in methodology is significant because the BBSOC method does not require a priori information about the solution structure.  In addition, the BBSOC method does not require user input and analysis to obtain an optimal solution. Also, as stated, the BBSOC method is implemented using a Legendre-Gauss-Radau (LGR) direct collocation method~\cite{Kameswaran2008,Garg2010,Garg2011a,Garg2011b,Patterson2015} where a convergence theory has been established under certain assumptions on smoothness and coercivity~\cite{Hager2016,Hager2017,Hager2018,Hager2019}.  Next, the BBSOC method produced more accurate solutions to the problem under consideration relative to the numerical solutions obtained in Ref.~\cite{Shen1999}.  Finally, as first observed in Ref.~\cite{Fleming2008}, three control torques are considered in this paper whereas only two control torques were considered in Ref.~\cite{Shen1999}.  The inclusion of a third control torque produces smaller optimal terminal times because the third control torque makes it possible to optimize the rotation of the body about the axis of symmetry and, thus, allowing the body to spin about its symmetry axis at a non-constant rate.

The remainder of the paper is organized as follows. Section~\ref{sect:probForm} presents the system model and problem formulation studied. Section~\ref{sect:mathAnalysis} presents the optimality conditions derived from Pontryagin's minimum principle along with a brief analysis of the singular controls. Section~\ref{sect:BBSOC} provides an overview of the numerical approach used to solve the problem, the BBSOC method, and the numerical approach developed in Ref.~\cite{Shen1999}. Section~\ref{sect:results} provides numerical results for four different maneuvers studied. The results are then summarized and further discussed in Section~\ref{sect:discussion}. Finally, Section~\ref{sect:conc} provides conclusions on this research.

\section{Problem Formulation}\label{sect:probForm}

In this section, the optimal control problem corresponding to the minimum-time reorientation of an axisymmetric rigid spacecraft is derived.  The formulation provided in this paper uses three control inputs and is similar to that given in~\cite{Shen1999}.  Section~\ref{sect:model} provides a derivation of the equations of motion for the spacecraft while Section~\ref{sect:ocp} provides the formulation of the optimal control problem.

\subsection{Equations of Motion for an Axisymmetric Rigid Spacecraft\label{sect:model}}

\begin{figure}[ht]
\begin{center}
	\includegraphics[scale=0.6]{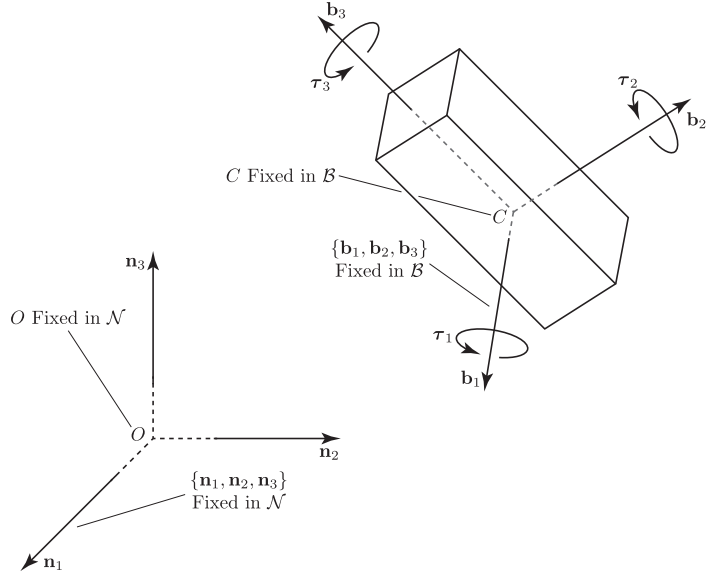}
	\caption{Axisymmetric rigid body with three control torques.}	
	\label{fig:rigidBody}
\end{center}
\end{figure}

Consider an axisymmetric rigid spacecraft as shown schematically in Fig.~\ref{fig:rigidBody}.  In order to model the motion of the spacecraft, the following bases are used.  First, let $\m{b}=\left\{\m{b}_1,\m{b}_2,\m{b}_3\right\}$ be a principal-axis basis that is fixed in the rigid body $\C{B}$ such that the direction $\m{b}_3$ lies along the axis of symmetry of the spacecraft, and let $C$ be the center of mass of $\C{B}$.  Then $\{C,\m{b}\}$ forms a principal-axis body-fixed coordinate system.  Next, let $\m{n}=\left\{\m{n}_1,\m{n}_2,\m{n}_3\right\}$ be a basis that is fixed in the inertial reference frame $\C{N}$, and let $O$ be fixed in $\C{N}$.  Then $\{O,\m{n}\}$ forms an inertially fixed coordinate system.  Furthermore, assume that $\C{B}$ is subject to actions by three control torques $\g{\tau}_1$, $\g{\tau}_2$, and $\g{\tau}_3$ about the directions $\m{b}_1$, $\m{b}_2$, and $\m{b}_3$, respectively.  Next, the angular velocity of the body relative to the inertial frame, $\lrsup{\C{N}}{\g{\omega}}{\C{B}}$, is defined as
\begin{equation}\label{angular-velocity}
  \lrsup{\C{N}\hspace{-2pt}}{\g{\omega}}{\C{B}} = \sum_{j=1}^3\omega_j\m{b}_j.
\end{equation}
Furthermore, let $(I_1,I_2,I_3)$ be the moments of inertia expressed in the body-fixed principal-axis coordinate system $\left\{C,\m{b}\right\}$, where $C$ is the center of mass of $\C{B}$. Finally $x_1$ and $x_2$, as introduced below, represent the relative position of the $\m{n}_3$ axis in the basis $\m{b}$ for a spin-stabilized spacecraft. The problem under consideration has the physical interpretation that for an observer in the body-fixed frame, the location of the $\m{n}_3$ axis is given by $x_1$ and $x_2$. The minimum-time problem aims to reorient the spacecraft between relative positions of the $\m{n}_3$ axis. Now the Eulerian dynamic equations of motion are introduced.

Applying a moment balance relative to the center of mass of the spacecraft gives
\begin{equation}\label{eq:EulersEqn}
  \m{M} = \m{I}^{\C{B}}\cdot\lrsup{\C{N}\hspace{-2pt}}{\g{\alpha}}{\C{B}} + \lrsup{\C{N}\hspace{-2pt}}{\g{\omega}}{\C{B}} \times \m{I}\lrsup{\C{N}}{\g{\omega}}{\C{B}},
\end{equation}
where $\m{I}^{\C{B}}$ is the moment of inertia tensor of $\C{B}$ relative to the center of mass of $\C{B}$, $\m{T}\cdot\m{c}$ is the operation that takes the product of a tensor $\m{T}$ with a geometric vector $\m{c}$, and 
\begin{equation}\label{angular-acceleration}
  \lrsup{\C{N}\hspace{-2pt}}{\g{\alpha}}{\C{B}} = \lsupddt{\C{N}\hspace{-3pt}}{}\left(\lrsup{\C{N}\hspace{-2pt}}{\g{\omega}}{\C{B}}\right) = \lsupddt{\C{B}\hspace{-3pt}}{}\left(\lrsup{\C{N}\hspace{-2pt}}{\g{\omega}}{\C{B}}\right) = \sum_{j=1}^3\dot{\omega}_j\m{b}_j, 
\end{equation}
is the angular acceleration of $\C{B}$ as viewed by an observer fixed in $\C{N}$.  Next, 
\begin{equation}\label{moment}
  \m{M} = \sum_{j=1}^3 \g{\tau}_j = \sum_{j=1}^3 \tau_j\m{b}_j
\end{equation}
where $\m{M}$ is the resultant moment applied to $\C{B}$ relative to $C$ and $\g{\tau}_j = \tau_j\m{b}_j,\; (j=1,2,3)$ are the pure torques applied to $\C{B}$ about the body-fixed directions $\m{b}_1$, $\m{b}_2$, and $\m{b}_3$, respectively.  
\begin{equation}\label{moment-of-inertia-tensor}
  \m{I}^{\C{B}} = \sum_{j=1}^3 I_j \m{b}_j\otimes \m{b}_j,
\end{equation}
is the moment of inertia tensor of the rigid body relative to $C$ (where $\m{a}\otimes\m{b}$ is the tensor product between two vectors $\m{a}$ and $\m{b}$).  Expressing all quantities in the body-fixed basis $\m{b}$ leads to the following three scalar differential equations for the rotational motion of the spacecraft:  
\begin{equation}\label{eq:SimpEulerEqn}
	\begin{array}{lcl}
		I_1\dot{\omega}_1 &=& (I_2-I_3)\omega_2\omega_3+\tau_1, \\
		I_2\dot{\omega}_2 &=& (I_3-I_1)\omega_3\omega_1+\tau_2, \\
		I_3\dot{\omega}_3 &=& (I_1-I_2)\omega_1\omega_2+\tau_3.
	\end{array}
\end{equation}
For an axisymmetric rigid body, $I_1=I_2$ when the axis of symmetry is chosen as the third axis which leads to the simplification of Eq.~\eqref{eq:SimpEulerEqn} as
\begin{equation}
	\begin{array}{lcl}
		\dot{\omega}_1 &=& \displaystyle \phantom{-} \frac{I_2-I_3}{I_1}\omega_{3}\omega_2+\frac{\tau_1}{I_1} \vspace{0.1cm}, \\ \vspace{0.1cm}
		\dot{\omega}_2 &=& \displaystyle -\frac{I_2-I_3}{I_1}\omega_{3}\omega_1+\frac{\tau_2}{I_2}, \\
		\dot{\omega}_3 &=& \displaystyle \phantom{-} \frac{\tau_3}{I_3}.
	\end{array}
\end{equation}
Furthermore, in Ref.~\cite{Shen1999}, $\omega_3$ is forced to remain constant throughout the maneuver, but in this study $\omega_3$ is included as a state and allowed to vary with the inclusion of the third control torque.
Finally, new control inputs can be defined as $u_j=\tau_j/I_j, \; (j = 1,2,3)$.  Then, defining
\begin{equation}
  a = (I_2-I_3)/I_1
\end{equation}
leads to the following three first-order differential equations:
\begin{equation}
  \begin{array}{lcl}
    \dot{\omega}_1 &=& \phantom{-}  a\omega_{3}\omega_2+u_1, \\
    \dot{\omega}_2 &=& -a\omega_{3}\omega_1+u_2, \\
    \dot{\omega}_3 &=& \phantom{-} u_3.
  \end{array}
\end{equation}
In order for the problem to be physically realizable it is necessary that $|a|\leq 1$.  

It is demonstrated in~\cite{Tsiotras1995}, that the parameterization of the attitude kinematics can be derived as follows. The orientation of the $\m{n}_3$ inertial axis as viewed by an observer in the body-fixed frame $\C{B}$ can be uniquely described by the variables $x_1$ and $x_2$ as
\begin{equation}
  \begin{array}{lcl}
    x_1 &=& \displaystyle \frac{\beta}{1+\gamma} ,\\
    x_2 &=& \displaystyle \frac{\sigma}{1+\gamma},
  \end{array}
\end{equation}
where $(\sigma,\beta,\gamma)$ are the direction cosines of $\m{n}_3$ with respect to $\m{b}$, i.e. $\m{n}_3 = \sigma \m{b}_1 + \beta \m{b}_2 + \gamma \m{b}_3$. Furthermore, it is shown in~\cite{Tsiotras1995} that $x_1$ and $x_2$ obey the following system of differential equations
\begin{equation}
  \begin{array}{lcl}
    \dot{x}_1 &=& \displaystyle \phantom{-} \omega_{3}x_2 + \omega_2x_1x_2 + \frac{\omega_1}{2}(1+x_1^2-x_2^2), \\
    \dot{x}_2 &=& \displaystyle -\omega_{3}x_1 + \omega_1x_1x_2 + \frac{\omega_2}{2}(1+x_2^2-x_1^2) .
  \end{array}
\end{equation}
Using the aforementioned parameterization, it is not possible to specify the absolute orientation of the rigid spacecraft in the inertial frame, nor is it possible to define the relative orientation of the spacecraft about the direction $\m{n}_3$ because a third parameter is required.  For spin-stabilized axisymmetric spacecraft the relative orientation of the body about the spin axis (or symmetry axis) is assumed to be irrelevant.  Consequently, the third parameter is not required and only the orientation of the symmetry axis is important.  Details on defining the third parameter for rigid body's using the parameterization considered are provided in Ref.~\cite{Tsiotras1995}.

\subsection{Optimal Control Problem\label{sect:ocp}}

Using the dynamic model developed in Section \ref{sect:model}, the optimal control problem that corresponds to the minimum-time reorientation of an axisymmetric rigid spacecraft can be posed as follows. Minimize the terminal time
\begin{equation}
  \C{J} = t_f
\end{equation}
subject to the dynamic constraint system 
\begin{equation}\label{eq:eom}
  \begin{array}{lcl}
    \dot{\omega}_1 &=& \phantom{-}  a\omega_{3}\omega_2+u_1, \\
    \dot{\omega}_2 &=& -a\omega_{3}\omega_1+u_2, \\
    \dot{\omega}_3 &=& \phantom{-}  u_3, \\
    \dot{x}_1 &=& \displaystyle \phantom{-}  \omega_{3}x_2 + \omega_2x_1x_2 + \frac{\omega_1}{2}(1+x_1^2-x_2^2) \vspace{0.05cm}, \\
    \dot{x}_2 &=& \displaystyle -\omega_{3}x_1 + \omega_1x_1x_2 + \frac{\omega_2}{2}(1+x_2^2-x_1^2), 
  \end{array}
\end{equation}
the control constraints
\begin{equation}
	u_{\min} \leq u_i \leq u_{\max}, \quad i = 1,2,3
\end{equation}
and the boundary conditions
\begin{equation}
  \begin{array}{lclclclcl}
    \g{\omega}(0) &=& [\omega_{10},\omega_{20},\omega_{30}]\T, \\
    \m{x}(0) &=& [x_{10}, x_{20}]\T ,\\
    \Phi[\g{\omega}(t_f),\m{x}(t_f)] &=& 0.
  \end{array}
\end{equation}
It is noted for this problem that the state, denoted $\m{Y}$, is given as
\begin{equation}
  \m{Y} = \left[ \begin{array}{c}
    \g{\omega} \\
    \m{x}
  \end{array}\right]
\end{equation}
and $\Phi:\R^5\rightarrow \R^k, \, k \leq 5$ is a differentiable vector function that defines the fixed and free terminal conditions of the state.

\section{Mathematical Analysis}\label{sect:mathAnalysis}

A mathematical analysis of the problem is now performed to analyze the solution structure using optimal control theory. First, the necessary conditions for optimality are obtained by applying Pontryagin's Minimum principle (PMP)~\cite{Kirk2004,Bryson1975,Athans2013}. First, the optimal control Hamiltonian is given as
\begin{equation}\label{eq:ham}
  \begin{split}
    \C{H} & = \lambda_{1} (a\omega_{3}\omega_2+u_1) + \lambda_{2} ( -a\omega_{3}\omega_1+u_2) + \lambda_3 u_3  
    + \lambda_{4} (\omega_{3}x_2 + \omega_2x_1x_2 + \frac{\omega_1}{2}(1+x_1^2-x_2^2)) \\ & + \lambda_{5} (-\omega_{3}x_1 + \omega_1x_1x_2 + \frac{\omega_2}{2}(1+x_2^2-x_1^2) ),
  \end{split}
\end{equation}
where $\g{\lambda}=(\lambda_1,\lambda_2,\lambda_{3},\lambda_{4},\lambda_5)$ is the costate vector and satisfies the differential equation 
\begin{equation}
  \dot{\g{\lambda}} = -\left[\frac{\partial \C{H}}{\partial\m{Y}} \right]\T.
\end{equation}
The differential equations for each component of the costate are then given, respectively, as
\begin{subeqnarray}\label{eq:costates}
  \slabel{eq:lam-1} \dot{\lambda}_{1} &=& \lambda_{2}a\omega_{3} - \lambda_{4}\frac{1+x_1^2-x_2^2}{2} -\lambda_{5}x_1x_2, \\
  \slabel{eq:lam-2} \dot{\lambda}_{2} &=& -\lambda_{1}a\omega_{3}  -\lambda_{4}x_1x_2 - \lambda_{5}\frac{1+x_2^2-x_1^2}{2}, \\
  \slabel{eq:lam-3} \dot{\lambda}_{3} &=& -\lambda_1 a\omega_2 + \lambda_2 a \omega_1 - \lambda_4 x_2 + \lambda_5 x_1,	\\		
  \slabel{eq:lam-4} \dot{\lambda}_{4} &=& -\lambda_{4} ( \omega_2 x_2 + \omega_1 x_1 ) +\lambda_{5} ( \omega_{3} - \omega_1 x_2 + \omega_2 x_1 ), \\
  \slabel{eq:lam-5} \dot{\lambda}_{5} &=& -\lambda_{4} ( \omega_{3} + \omega_2 x_1 - \omega_1 x_2 ) - \lambda_{5} ( \omega_1 x_1 + \omega_2 x_2 ) ,
\end{subeqnarray}
where the terminal costate $\g{\lambda}(t_f)$ is defined by the transversality condition
\begin{equation}\label{eq:transversalityCond}
  \g{\lambda}\T(t_f) = \g{\nu}\T\frac{\partial\Phi}{\partial \m{Y}(t_f)}
\end{equation}
and $\g{\nu} \in \R^k$ is a constraint multiplier.  Next, the transversality condition associated with the terminal time is 
\begin{equation}
  \C{H}(t_f) = -1
\end{equation}
Furthermore, the three components of the optimal control are obtained from the strong form of Pontryagin's Minimum Principle as
\begin{equation}\label{eq:PMPstrong}
  \m{u}_j^*(t) = \min_{u_j \in [u_{\min},u_{\max}]} \C{H}(\m{Y}(t),u_i), \quad \forall t \geq 0, \quad  (j =1,2,3).  
\end{equation}
Because the control components appear linearly in the Hamiltonian, each component of the optimal control is given as
\begin{equation}\label{eq:PMPweak}
  u_j^*
  =\left\{
    \begin{array}{lcl}
      u_{\min} & , & g_j>0, \\
      u_{j,\textrm{s}} & , & g_j=0, \\
      u_{\max} & , & g_j<0,
    \end{array}
    , \quad \quad  (j = 1,2,3),
    \right.
\end{equation}
where
\begin{equation}\label{eq:switching-functions}
  g_j = \lambda_j, \quad (j=1,2,3)
\end{equation}
are the switching functions corresponding to the components $u_j,\;(j=1,2,3)$ of the control, respectively.

\subsection{Singular Control Analysis\label{sect:singArcAnalysis}}

For the following discussion, it is assumed that $a\neq 0$ and $\omega_{3}\neq 0,\, \forall t\in[t_0,t_f]$. The switching functions defined in Eq.~\eqref{eq:switching-functions} determine the possible occurrences of a singular arc in the solution. The control $u_j^*$ is singular whenever $g_j = 0$ during a nonzero interval $[t_1,t_2] \subset [t_0,t_f]$. The singular control is obtained implicitly from the switching function. Specifically, $g_j$ is differentiated repeatedly until the control $u_j$ explicitly appears~\cite{Schattler2012}. Thus $u_j^*$ can be solved for by
\begin{equation}
	\frac{d^{(2r_j)}}{dt^{(2r_j)}}\left(g_j\right)= 0, \quad (r_j = 0,1,2,\ldots), \quad j= 1,2,3
\end{equation}
where $2r_j$ is the minimum number of differentiation's of $g_j$ required to obtain the corresponding control $u_{j,\textrm{s}}$.
For $u_j$ to be optimal over a singular arc, the number of differentiation's $2r_j$ must be even~\cite{Bryson1975,Schattler2012}.
Furthermore, the generalized Legendre-Clebsch condition~\cite{Bryson1975,Kelley1967,Kopp1965}
\begin{equation}\label{eq:gen-Legendre-Clebsch}
	(-1)^r \frac{\partial}{\partial u_j}\left[ \frac{d^{2r}}{dt^{(2r)}}g_j \right] \geq 0, \quad (r_j = 0,1,2,\ldots), \quad j = 1,2,3
\end{equation}
must hold over the duration of a singular arc $[t_1,t_2]$.

An in depth analysis of the singular controls in the minimum-time reorientation problem was performed in~\cite{Shen1999}. The analysis in Ref.~\cite{Shen1999} is extended to consider three torque controls as opposed to two. When considering the possibility of all three controls being singular, the analysis remains the same where all three controls being singular is not optimal. As it is shown in Ref.~\cite{Shen1999}, the only optimal case to consider is the scenario where one control component is singular and the other two control components are bang-bang. Furthermore, the authors were not able to find a set of boundary conditions that resulted in the third control component being singular. Thus, only the case where either the first or second control component is singular is considered in this work. A brief review of this scenario is provided below including the extension to three control torques. Special cases of the minimum-time reorientation problem will also be analyzed later in Section~\ref{sect:specialCases}.

\subsection{Either \texorpdfstring{$u_1$}{TEXT} or \texorpdfstring{$u_2$}{TEXT} is Singular While Other Control Components are Bang-Bang}

Consider the case where one of the first two control components is singular while the other two control components are bang-bang.  In any of these cases, the singular control is obtained by differentiating the switching function corresponding to the singular control component with respect to time until the control appears.  As it turns out, the switching function of interest needs to be differentiated with respect to time {\em four} times until the corresponding control component appears, see Appendix~\ref{sect:appGeneral} for the derivation. The singular arc is of second order given that four time differentiation are needed. The singular control components $u_{1,\textrm{s}}$ or $u_{2,\textrm{s}}$ are then given as follows: 
\begin{equation}\label{eq:singularLawU1}
	u_{1,\textrm{s}} = -\frac{\dot{\lambda}_2[ 2a\omega_{3}^3(a+1)^2 - 2a\omega_{3}\omega_1^2 + 2a\omega_{3} \omega_2^2 + 2\omega_1u_2 ] - \omega_2\lambda_2(4a^2\omega_{3}^2\omega_1 - 3a\omega_{3}u_2 ) }{\dot{\lambda}_2 \omega_2 - \omega_{3}(1+2a)(\lambda_3x_2 - \lambda_4x_1)  }
\end{equation}
\begin{equation}\label{eq:singularLawU2}
	u_{2,\textrm{s}} = -\frac{\dot{\lambda}_1[ 2a\omega_{3}^3(a+1)^2 - 2a\omega_{3}\omega_1^2 + 2a\omega_{3} \omega_2^2 + 2\omega_2u_1 ] - \omega_1\lambda_1(4a^2\omega_{3}^2\omega_2 + 3a\omega_{3}u_1 ) }{\dot{\lambda}_1 \omega_1 - \omega_{3}(1+2a)(\lambda_3x_2 - \lambda_4x_1)  }
\end{equation}
The remaining bang-bang control components are determined according to Eq.~\eqref{eq:PMPweak}.
Furthermore, the generalized Legendre-Clebsch condition
\begin{equation}
  \frac{\partial}{\partial u_j}\left[ \frac{d^{4}}{dt^{4}}g_j \right] \geq 0, \quad j = 1,2
\end{equation}
implies that the denominator of Eqs.~\eqref{eq:singularLawU1}--\eqref{eq:singularLawU2} must be nonnegative in order for the singular control to be optimal.  

\subsection{Special Cases: Either \texorpdfstring{$\omega_3$}{TEXT} or \texorpdfstring{$a$}{TEXT} Is Zero \label{sect:specialCases}}

Two special cases arise when either $a$ or $\omega_3$ is zero.  When $\omega_{3}= 0$, the equations of motion given in Eq.~\eqref{eq:eom} are simplified greatly.  When $a=0$ and $\omega_3$ is constant, the equations of motion given in Eq.~\eqref{eq:eom} simplify to a system that represents an inertially symmetric rigid body.  The cases $\omega_3= 0$ and $a=0$ are analyzed in Sections.~\ref{sect:nonspinningCase} and~\ref{sect:inertialCase}, respectively.

\subsubsection{\texorpdfstring{$\omega_3=0$}{TEXT}: Nonspinning Axisymmetric Rigid Body\label{sect:nonspinningCase}}

When $\omega_{3}=0,\, \forall t\in[t_0,t_f]$, Eq.~\eqref{eq:eom} simplifies to
\begin{equation}\label{eq:FS-eom}
  \begin{array}{lcl}
    \dot{\omega}_1 &=& u_1 \\
    \dot{\omega}_2 &=& u_2 \\
    \dot{\omega}_3 &=& 0 \\
    \dot{x}_1 &=& \omega_2x_1x_2 + \frac{\omega_1}{2}(1+x_1^2-x_2^2) \\
    \dot{x}_2 &=& \omega_1x_1x_2 + \frac{\omega_2}{2}(1+x_2^2-x_1^2) 
  \end{array}
\end{equation}
and the resulting analysis for $u_1$ singular and $u_2$ bang-bang in the section above is simplified significantly. Note that $u_3$ is no longer considered since the angular velocity about the third-axis is zero. Specifically, the second-order singular optimal control is now defined as
\begin{equation}\label{eq:FSA}
	u_1^*= 0
\end{equation}
and the bang-bang control $u_2$ is defined as
\begin{equation}\label{eq:bangbangu2}
  u_2^*
    =\left\{
    \begin{array}{lcl}
      u_{2,\min} & , & \lambda_2>0, \\
      u_{2,\max} & , & \lambda_2<0,
    \end{array}
    \right.
\end{equation}
See Appendix~\ref{sect:appSpecial} for the derivation of Eq.~\eqref{eq:FSA}.

\subsubsection{\texorpdfstring{$a=0$}{TEXT}:  Inertially Symmetric Rigid Body\label{sect:inertialCase}}

When $a=0$ and $\omega_3$ is constant, Eq.~\eqref{eq:eom} simplifies to 
\begin{equation}\label{eq:IS-eom}
	\begin{array}{lcl}
		\dot{\omega}_1 &=& u_1 \\
		\dot{\omega}_2 &=& u_2 \\
		\dot{\omega}_3 &=& 0 \\
		\dot{x}_1 &=& \omega_{3}x_2 + \omega_2x_1x_2 + \frac{\omega_1}{2}(1+x_1^2-x_2^2) \\
		\dot{x}_2 &=& -\omega_{3}x_1 + \omega_1x_1x_2 + \frac{\omega_2}{2}(1+x_2^2-x_1^2) 
	\end{array}
\end{equation}
and, again, the resulting analysis for $u_1$ singular and $u_2$ bang-bang in the section above is simplified significantly. For this scenario, $\omega_3$ spins at a constant rate and $\dot{\omega}_3=0$. Thus $u_3$ is no longer considered. In this case it is shown in Ref.~\cite{Shen1999} that all higher-order derivatives of $g_1$ are identically zero resulting in an optimal control $u_1$ that lies on an infinite-order singular arc. Thus, any control within the limits $[u_{1,\min},u_{1,\max}]$ that satisfy the boundary conditions qualifies as an optimal control. This result leads to multiple solutions for the infinite-order singular control that will minimize the final time of the problem. As in the previous case, the bang-bang control $u_2$ is determined by Eq.~\eqref{eq:bangbangu2}.

\section{Numerical Approach: BBSOC Method\label{sect:BBSOC}}

The minimum-time axisymmetric spacecraft reorientation problem studied in this paper is solved using the bang-bang and singular optimal control (BBSOC) method developed in Ref.~\cite{Pager2022}. The BBSOC method can identify the existence of singular arcs in a control solution and then perform an iterative regularization procedure on an interval where a singular control is detected. The BBSOC method can also identify the switching structure of a nonsmooth control solution and optimize values of the switch times. Specifically, the structure identification procedure of the BBSOC method analyzes an initial solution to an optimal control problem for any discontinuities, bang-bang arcs, and/or singular arcs. After the structure has been identified, the BBSOC method partitions the initial mesh into domains representing each identified arc.  Finally, corresponding constraints are enforced in each domain, respectively, while a regularization procedure is applied in any domains that have been identified as singular. The method is algorithmic in nature and requires no user-input.  A flow chart for the BBSOC method is provided in Fig.~\ref{fig:bbsoc-flowchart}. In addition, the BBSOC method is implemented using a multiple-domain formulation of the $hp$-adaptive Legendre-Gauss-Radau (LGR) collocation method~\cite{Garg2010,Garg2011a,Garg2011b,Kameswaran2008,Patterson2015}, which allows the use of variable mesh refinement to exploit knowledge of the identified control structure. More details on the BBSOC method and the multiple-domain LGR collocation can be found in Refs.~\cite{Pager2022,Pager2022conf}.

\subsection{Comparison with Ref.~\cite{Shen1999}\label{sect:methodComp}}

An important contribution of this paper is comparing the approach used in this paper with the approach used in Ref.~\cite{Shen1999}.  Specifically, the BBSOC method has been developed for solving nonsmooth and singular optimal control problems. In this paper, the performance of the BBSOC method is compared to the approach used in Ref.~\cite{Shen1999}.  Specifically, the method of Ref.~\cite{Shen1999} uses an indirect shooting method to solve the optimal control problem.  Given that a very good initial guesses is often required in order to successfully employ an indirect method, Ref.~\cite{Shen1999} first uses a direct method -- EZopt -- to obtain an initial guess for the indirect method -- BNDSCO -- used in Ref.~\cite{Shen1999}.  In particular, the Lagrange multipliers obtained from EZopt are used to obtain initial guesses for the costates using COSCAL~\cite{Shen1999}.  Then, an initial guess of the states, costates, controls along with the switching structure is provided to BNDSCO.

\begin{figure}[ht]
\centering
\vspace*{0.25cm}
\subfloat[\label{fig:bbsoc-flowchart}]{\includegraphics[scale=0.47]{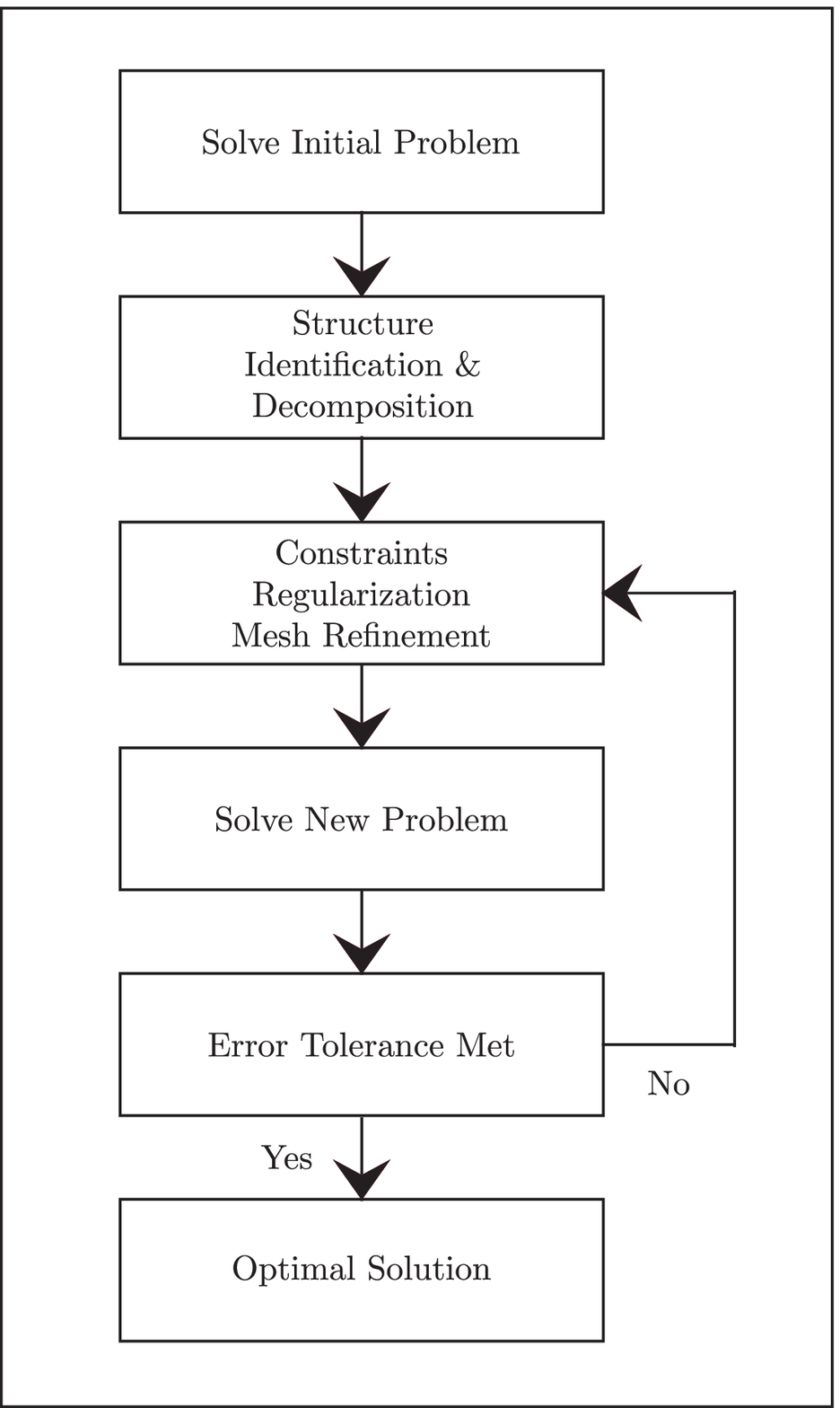}}~~\subfloat[\label{fig:shen-flowchart}]{\includegraphics[scale=0.47]{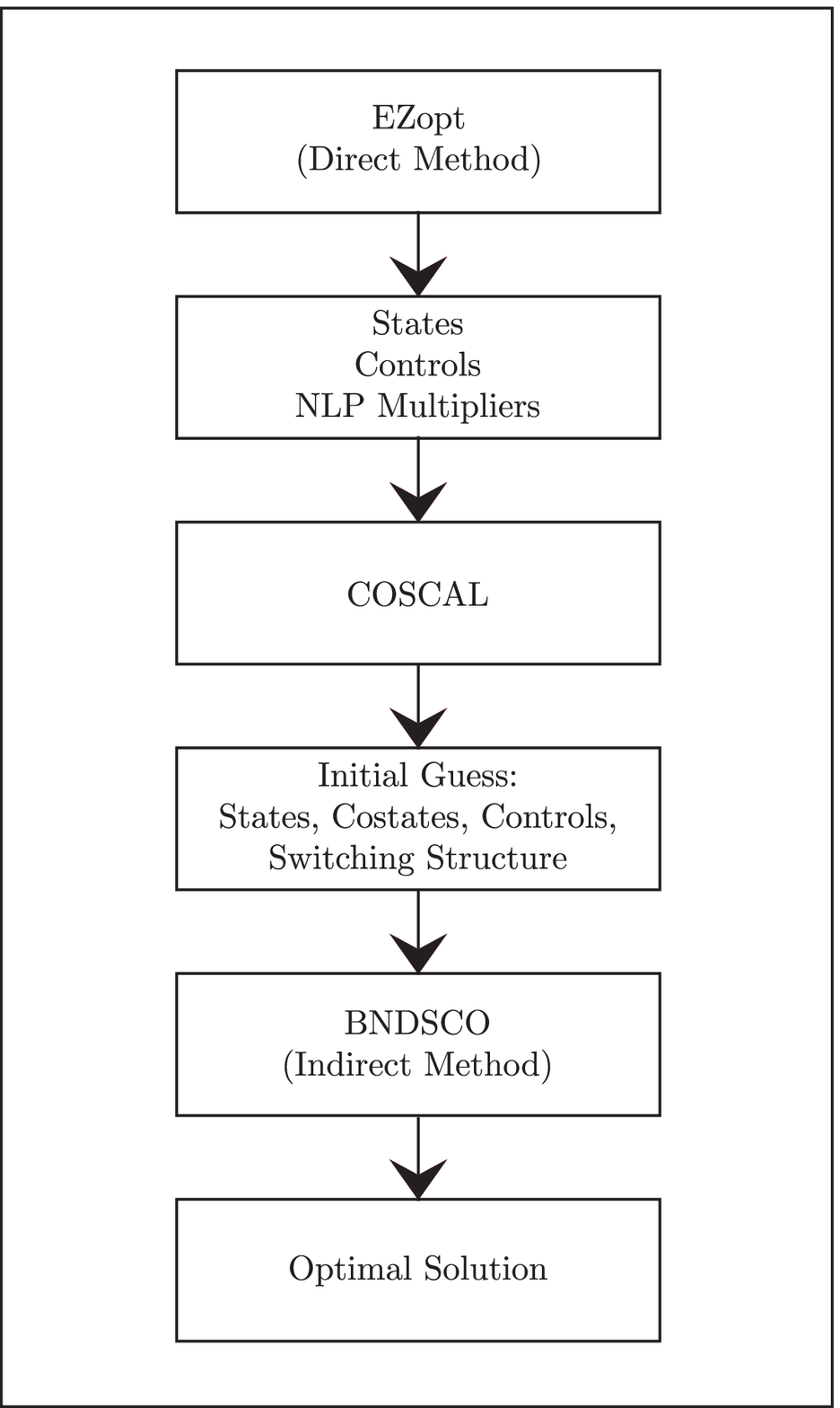}} \\
\caption{Flow chart of the (a) BBSOC method and (b) computational scheme used in Ref.~\cite{Shen1999}.} 
\label{fig:numApproach}
\end{figure}

A flow chart of the numerical approach of Ref.~\cite{Shen1999} is shown in Fig.~\ref{fig:shen-flowchart}. Given that an indirect shooting method is used to solve the optimal control problem, it is required that the higher-order singular arc optimality conditions be derived and that the switching structure of the control be known a priori.  For a general problem, however, deriving the higher-order optimality conditions may be difficult and the switching structure will not be known.  
On the other hand, the BBSOC method used to generate the results in this paper neither requires that the higher-order singular arc optimality conditions be derived nor does it require a priori knowledge of the switching structure of the optimal control .  Instead, the BBSOC method identifies algorithmically the structure of the optimal control and solves for the singular control through the aforementioned regularization procedure (which is itself algorithmic and requires no user intervention).  Figure~\ref{fig:numApproach} provides a flow chart that compares the BBSOC method to the method of Ref.~\cite{Shen1999}.  

\section{Results}\label{sect:results}

All results presented in this section were obtained using the bang-bang and singular optimal control (BBSOC) method~\cite{Pager2022} implemented in MATLAB\textsuperscript{\textregistered} with the NLP problem solver IPOPT~\cite{Biegler2008}.  It is noted that IPOPT was employed in full-Newton (second-derivative) mode with default NLP solver tolerance of $\epsilon_{NLP}=10^{-7}$. Any necessary mesh refinement was performed using the $ph$ mesh refinement described in Ref.~\cite{Patterson2015} with a mesh refinement accuracy tolerance of $\epsilon_{\textrm{mesh}}=10^{-5}$.  Furthermore, the lower and upper limits on the number of Legendre-Gauss-Radau points allowed in the method of Ref.~\cite{Patterson2015}  was $3$ and $12$, respectively.  The BBSOC method was initialized using a mesh that consisted of $20$ uniformly spaced mesh intervals and $3$ collocation points per mesh interval. A straight line guess is used for variables with boundary conditions at both endpoints and a constant guess is used for variables with boundary conditions at only one endpoint.  All first and second derivatives required by IPOPT were provided using the algorithmic differentiation software {\em ADiGator}~\cite{Weinstein2017}. Finally, all computations were performed on a 2.9 GHz 6-Core Intel Core i9 MacBook Pro running Mac OS Big Sur Version 11.6 with 32 GB 2400 MHz DDR4 of RAM, using MATLAB version R2019b (build 9.7.0.1190202) and all computation (CPU) times are in reference to this aforementioned machine.

Four different maneuvers of the minimum-time axisymmetric rigid spacecraft reorientation problem are presented and discussed below. The four cases include a rest-to-rest (RTR) maneuver where the resulting control is bang-bang, a non-rest-to-rest (NRTR) maneuver where the resulting control is bang-bang, a NRTR maneuver where the resulting control exhibits a finite-order singular arc, and a rest-to-non-rest (RTNR) maneuver where the resulting control exhibits an infinite-order singular arc. Each case is solved using the BBSOC method and $\mathbb{GPOPS-II}$~\cite{Patterson2014} (referred to as $hp$--LGR in discussion), and also compared with the results presented in Ref.~\cite{Shen1999} where the numerical approach discussed in Section~\ref{sect:methodComp} is implemented.

\subsection{Rest-to-Rest Maneuver: Bang-Bang Control\label{sect:bangSoln-RTR}}

The rest-to-rest (RTR) maneuver assumes both $a$ and $\omega_{3}$ are non-zero and $a=0.5$. The boundary conditions for this case are given as
\begin{displaymath}
  \begin{array}{lclclcl}
    \omega_{10} & = & 0 &,& \omega_{1f} & = & 0, \\
    \omega_{20} & = & 0 &,& \omega_{2f} & = & 0, \\
    \omega_{30} & = & -0.5 &,& \omega_{3f} & = & -0.5, \\
    x_{10} & = & 1.5 &,& x_{1f} & = & 0,  \\
    x_{20} & = & -0.5 &,&  x_{2f} & = & 0.  
  \end{array}
\end{displaymath}
For this maneuver, most RTR maneuvers, and non-rest-to-rest (NRTR) maneuvers where the body is spinning and axisymmetric, all three controls are bang-bang~\cite{Shen1999}.  This behavior is evident in the results of Section~\ref{sect:bangSoln-NRTR} where a NRTR maneuver is analyzed. For the NRTR and RTR maneuvers, the goal is to rotate the body until the direction $\m{b}_3$ is aligned with the direction $\m{n}_3$ representing a rest condition. The three components of the optimal control were found to be bang-bang such that $u_1$ contains a single switch while $u_2$ and $u_3$ contain two switches. The minimum time to complete the maneuver was approximately $2.5126$ seconds.

\begin{figure}[ht]
\centering
\vspace*{0.25cm}
\subfloat[\label{fig:bangRTR-U1}]{\includegraphics[scale=0.3]{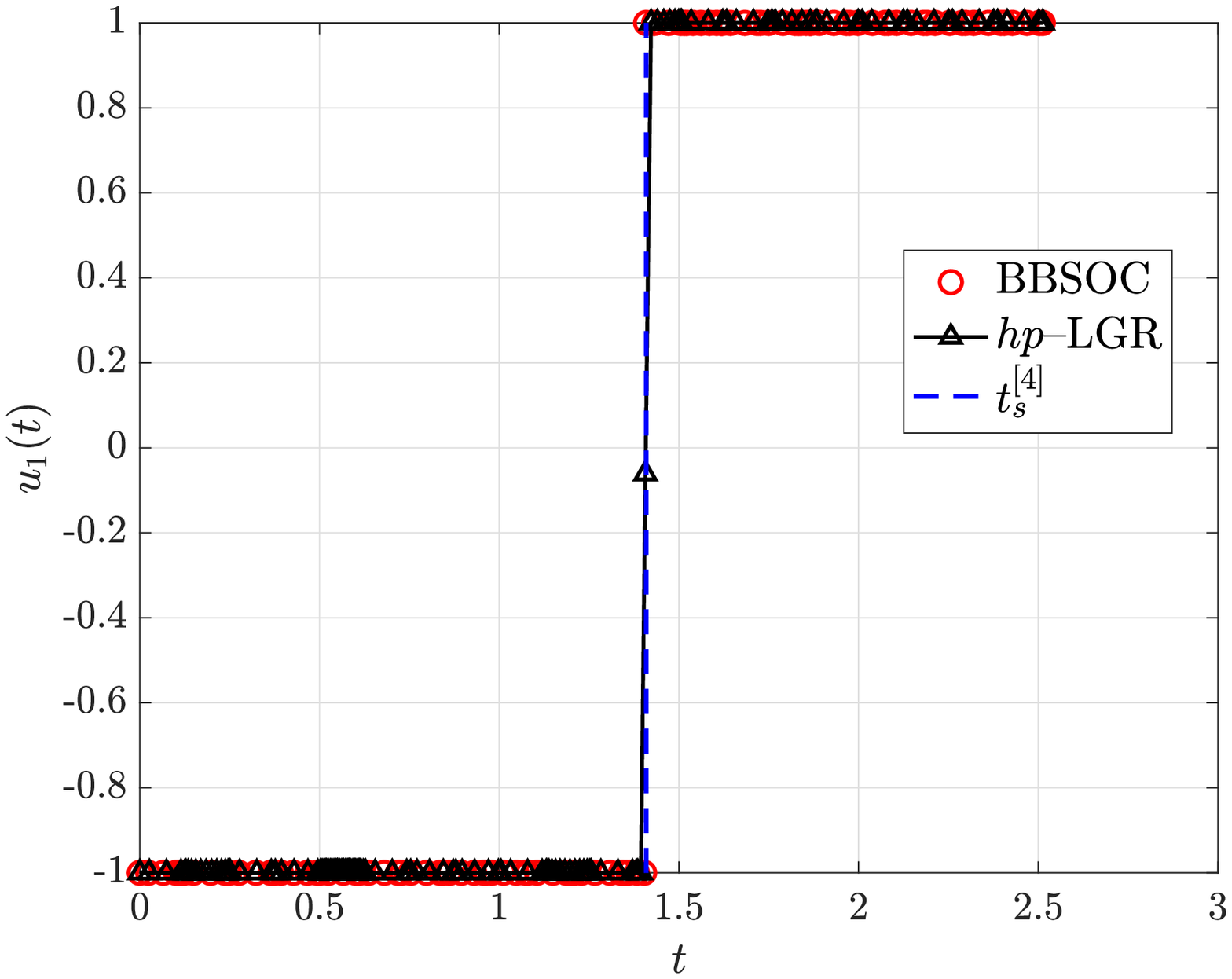}}~~\subfloat[\label{fig:bangRTR-U2}]{\includegraphics[scale=0.3]{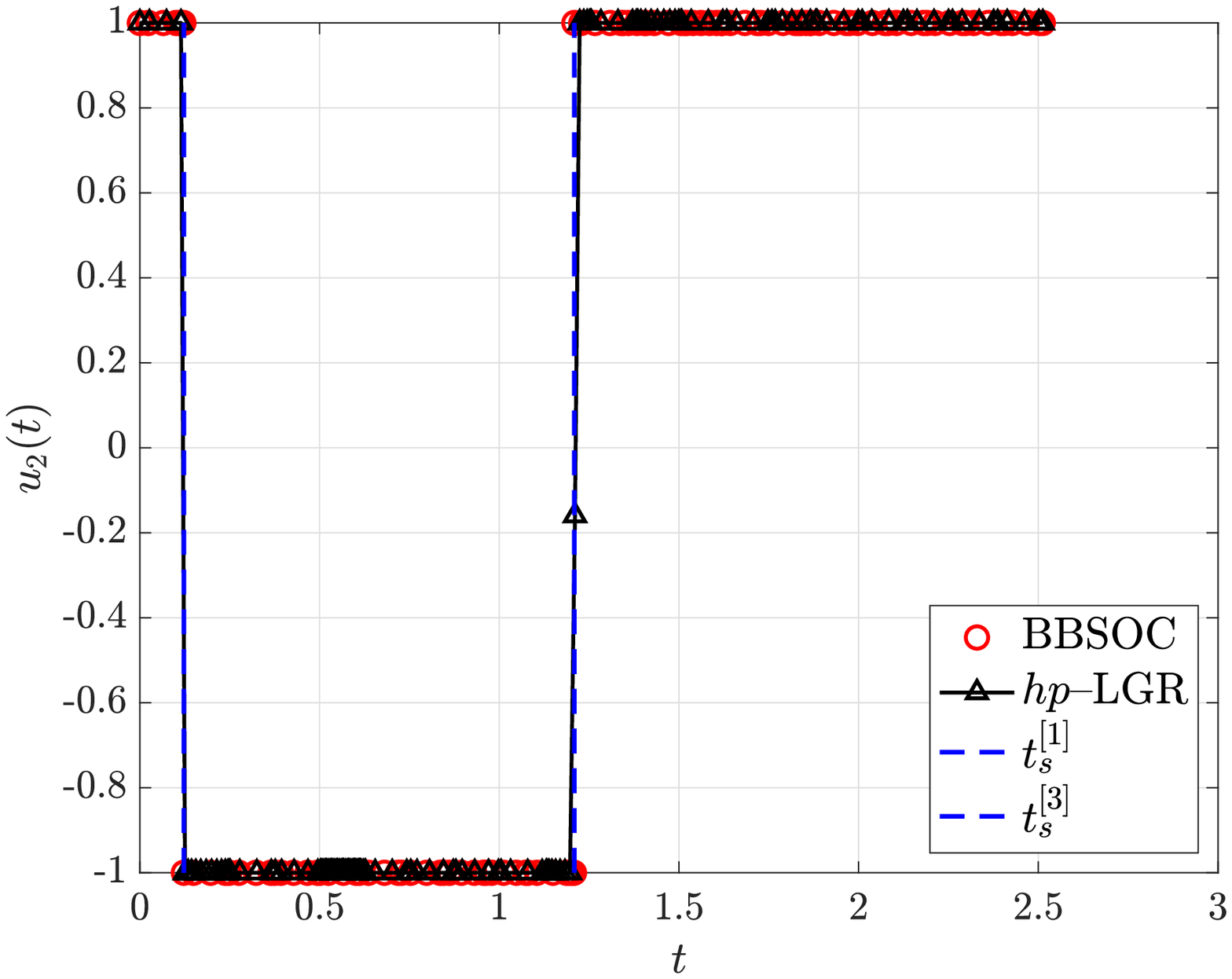}}~~\subfloat[\label{fig:bangRTR-U3}]{\includegraphics[scale=0.3]{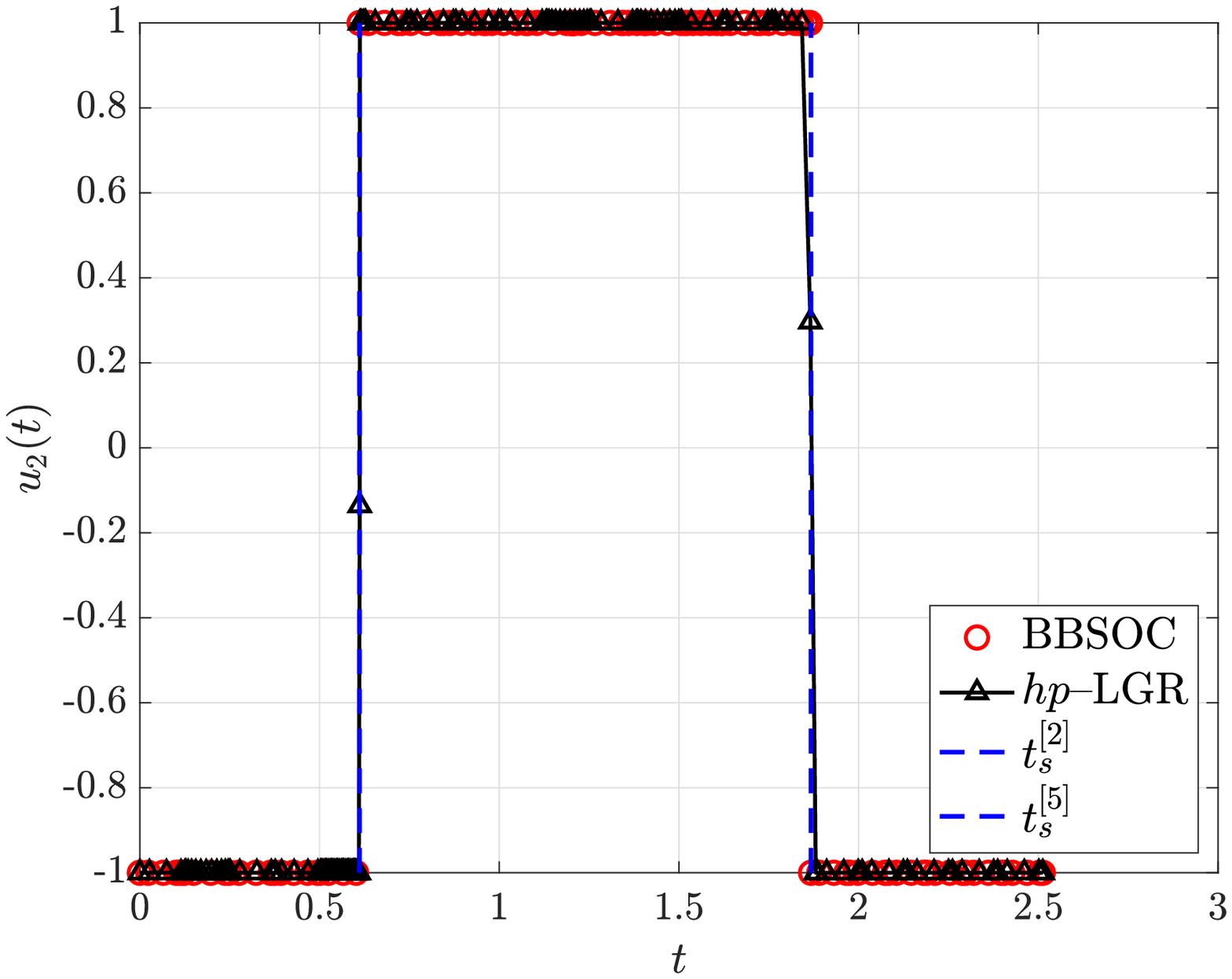}}  \\ 
\caption{Control component solutions for the spinning RTR maneuver. \label{fig:bangRTR1}}
\end{figure}

\begin{figure}[ht]
\centering
\vspace*{0.25cm}
\subfloat[\label{fig:bangRTR-PHI}]{\includegraphics[scale=0.3]{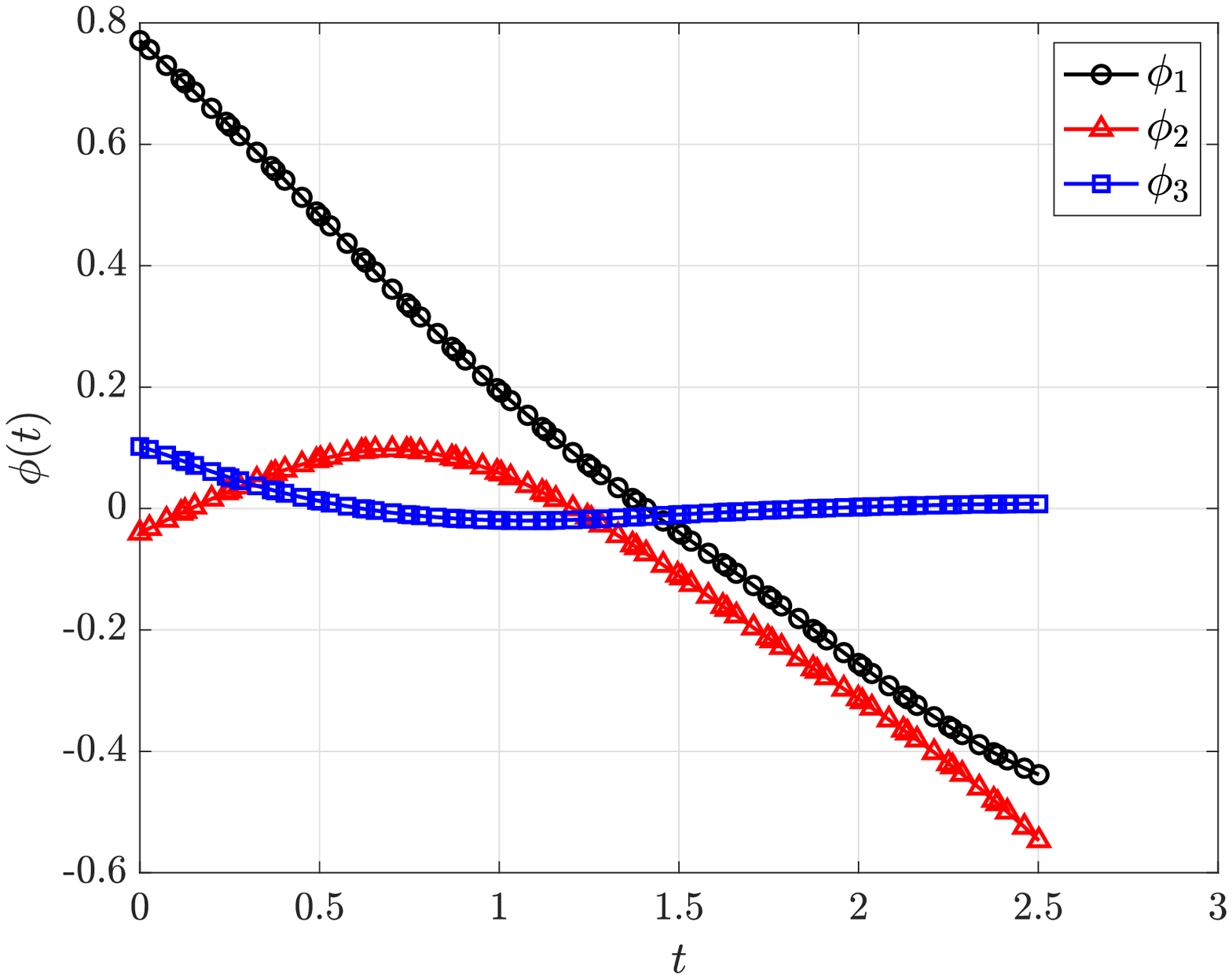}}~~ 
\subfloat[\label{fig:bangRTR-omega}]{\includegraphics[scale=0.3]{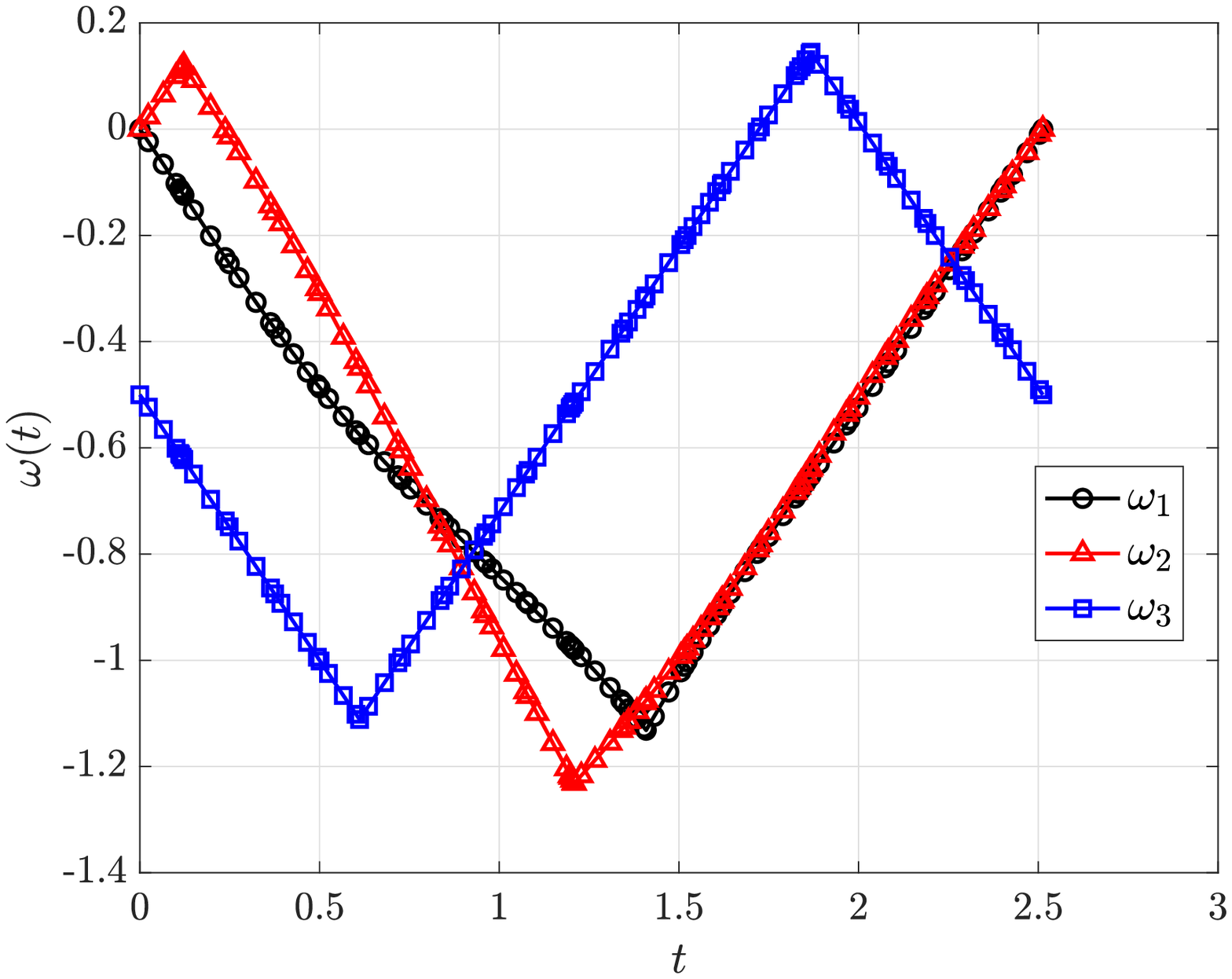}} \\
\caption{Switching functions and state solutions for the spinning RTR maneuver. \label{fig:bangRTR2}}
\end{figure}

Figures~\ref{fig:bangRTR1} and~\ref{fig:bangRTR2} shows the three control solutions along with their corresponding switching functions and the time history of the angular velocities. It is observed in Fig.~\ref{fig:bangRTR1}  that the BBSOC method identified the five switch times and constrained the control to the appropriate limit.  In contrast, the solution obtained using the $hp$-LGR method contains collocation points located between the lower and upper limits on the control in the neighborhood of the switch times. Table~\ref{tab:bangRTR} provides a further comparison between the BBSOC and $hp$-LGR methods by showing the performance of each method and by comparing the values of the switch points and the minimum terminal time.  Table~\ref{tab:bangRTR} also compares the results of Ref.~\cite{Shen1999}, where it is seen that the BBSOC method not only obtains a lower optimal terminal time, but also obtains the solution in a more computationally efficient manner than the method of Ref.~\cite{Shen1999}. It is noted that Ref.~\cite{Shen1999} claimed it took less than two seconds for BNDSCO to converge.  Note, however, that it took EZopt two minutes to converge to the initial guess required as the initial guess for BNDSCO.  As a result, it is seen that one of the major drawbacks of using an indirect method beyond requiring the derivation of the singular optimality conditions, is an increased computation time relative to a direct collocation method such as the BBSOC method.  Furthermore, by allowing $\m{b}_3$ to spin at a non-constant rate, a smaller terminal time is achieved by the three-torque control formulation used in this paper as compared to the two-torque formulation in Ref.~\cite{Shen1999}.  Specifically, Table~\ref{tab:bangRTR} shows a $3.73\%$ reduction in the terminal time obtained by the BBSOC method compared with the method of Ref.~\cite{Shen1999}.  

\begin{table}[ht]
\caption{\label{tab:bangRTR} Comparison of computational results for the RTR maneuver.}
\centering
\begin{tabular}{ c " c c c c c c c c c c c } 
\hline \hline
                &  $t_{s}^{[1]}$  & $t_{s}^{[2]}$   & $t_{s}^{[3]}$ & $t_{s}^{[4]}$ & $t_{s}^{[5]}$ & $t_f$ & CPU~$[s]$   \\ \thickhline 
BBSOC  & $0.1224$  & $0.6114$ & $1.2091$  & $1.4088$ & $1.8676$  & $2.5126$  &  $2.09$   \\ 
$hp$-LGR  & $0.1145$  & $0.6118$  & $1.2108$  & $1.4069$ & $1.8658$  & $2.5126$  & $4.41$ \\
Ref.~\cite{Shen1999}   & $-$  & $-$  & $-$  & $-$ & $-$  & $2.61$  & $122$  \\
\hline \hline
\end{tabular}
\end{table}

\subsection{Non-Rest-to-Rest Maneuver:  Bang-Bang Control}\label{sect:bangSoln-NRTR}

The non-rest-to-rest (NRTR) maneuver assumes both $a$ and $\omega_{3}$ are non-zero and $a=0.5$. The boundary conditions for this case are given as
\begin{displaymath}
  \begin{array}{lclclcl}
    \omega_{10} & = & -0.45 &,& \omega_{1f} & = & 0, \\
    \omega_{20} & = & -1.1 &,& \omega_{2f} & = & 0, \\
    \omega_{30} & = & -0.5 &,& \omega_{3f} & = & -0.5, \\
    x_{10} & = & 1.5 &,& x_{1f} & = & 0,  \\
    x_{20} & = & -0.5 &,&  x_{2f} & = & 0, 
  \end{array}
\end{displaymath}
and, as was the case in Section \ref{sect:bangSoln-RTR}, the optimal controls are bang-bang. The results obtained for this maneuver are similar to those presented in Section~\ref{sect:bangSoln-RTR}. Again, the time-optimal controls were found to be bang-bang with $u_1$ having one switch and $u_2$ and $u_3$ having two switches. The minimum time to complete the maneuver was $2.2445$ seconds.

\begin{figure}[ht]
\centering
\vspace*{0.25cm}
\subfloat[\label{fig:bangNRTR-U1}]{\includegraphics[scale=0.3]{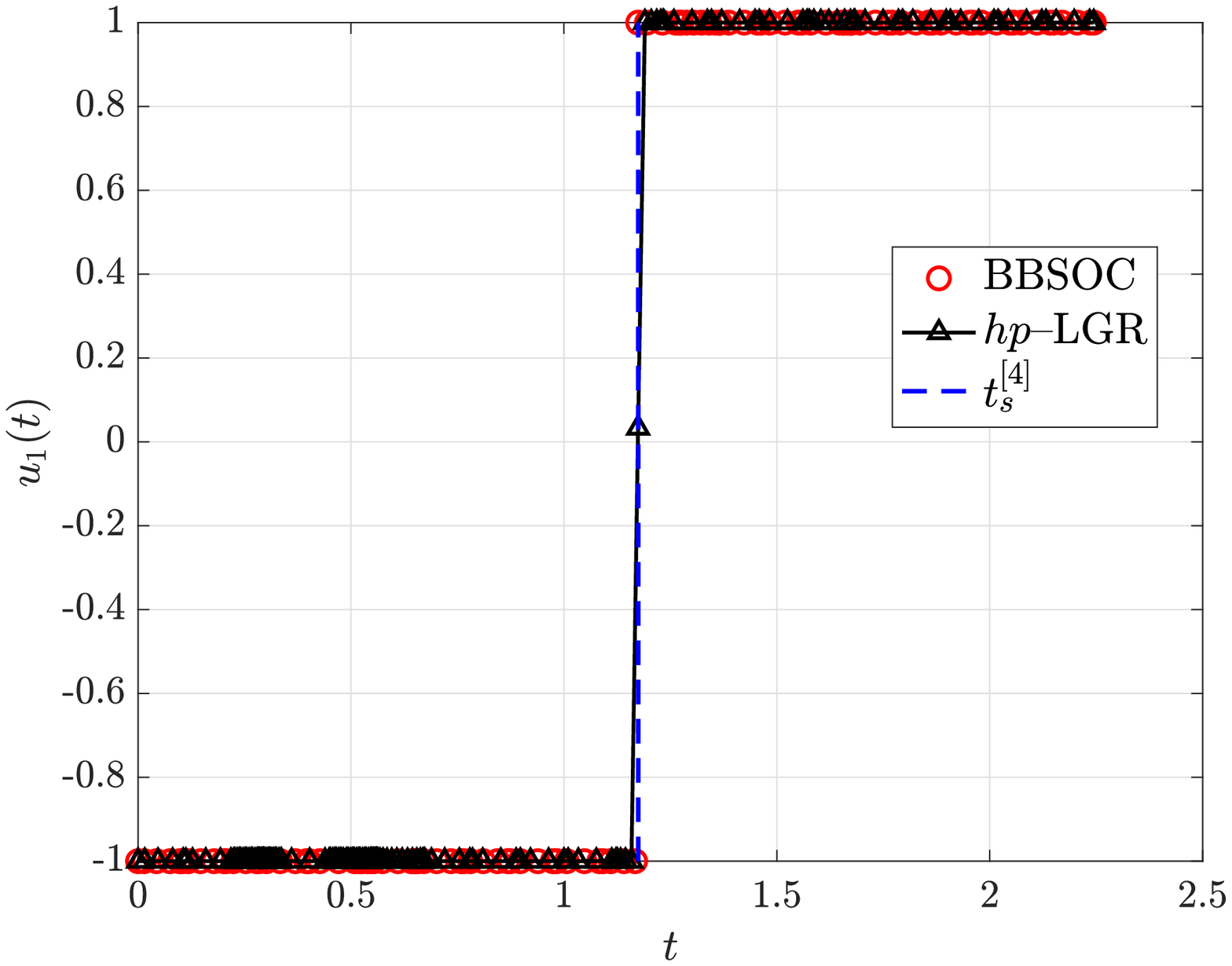}}~~\subfloat[\label{fig:bangNRTR-U2}]{\includegraphics[scale=0.3]{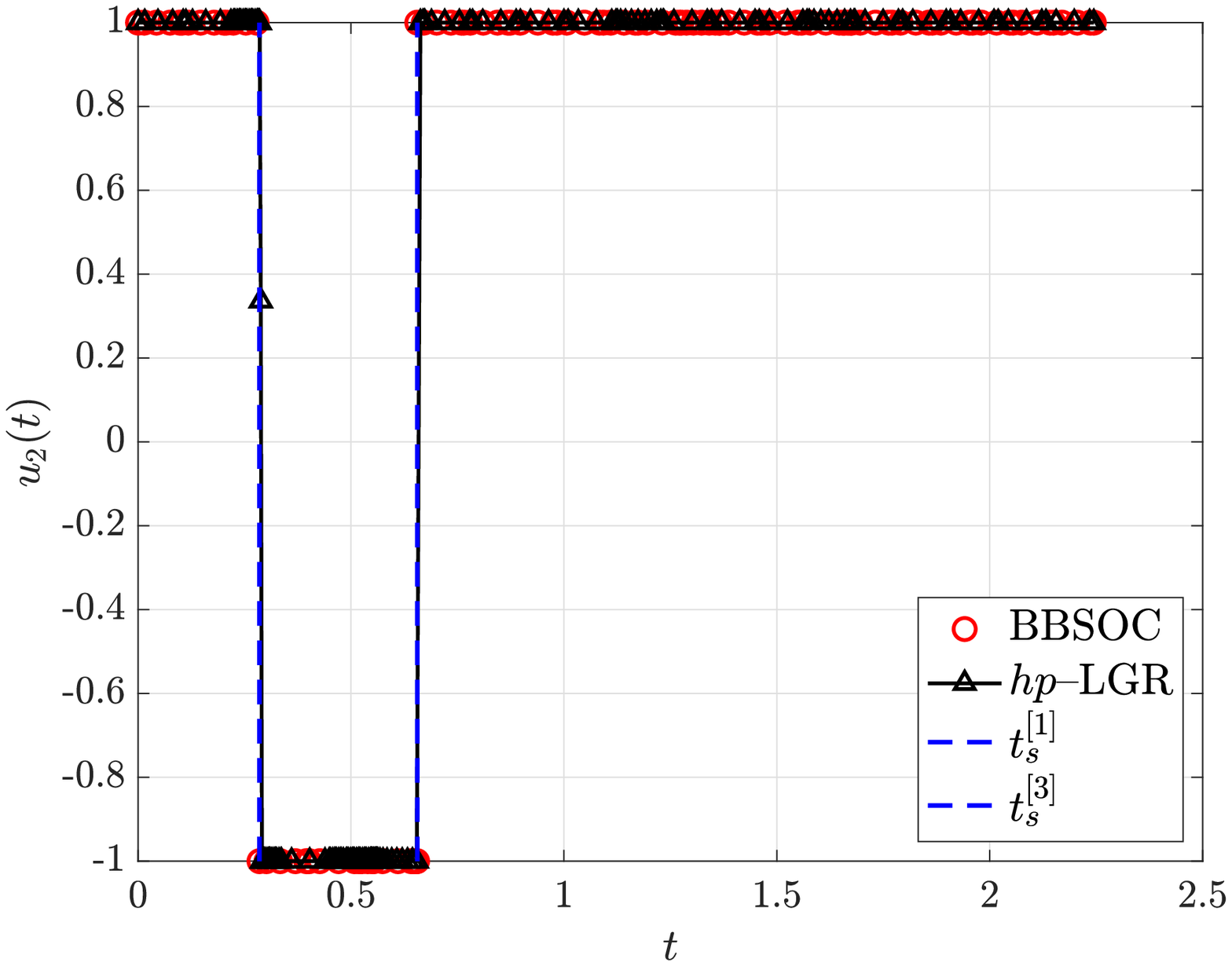}}~~\subfloat[\label{fig:bangNRTR-U3}]{\includegraphics[scale=0.3]{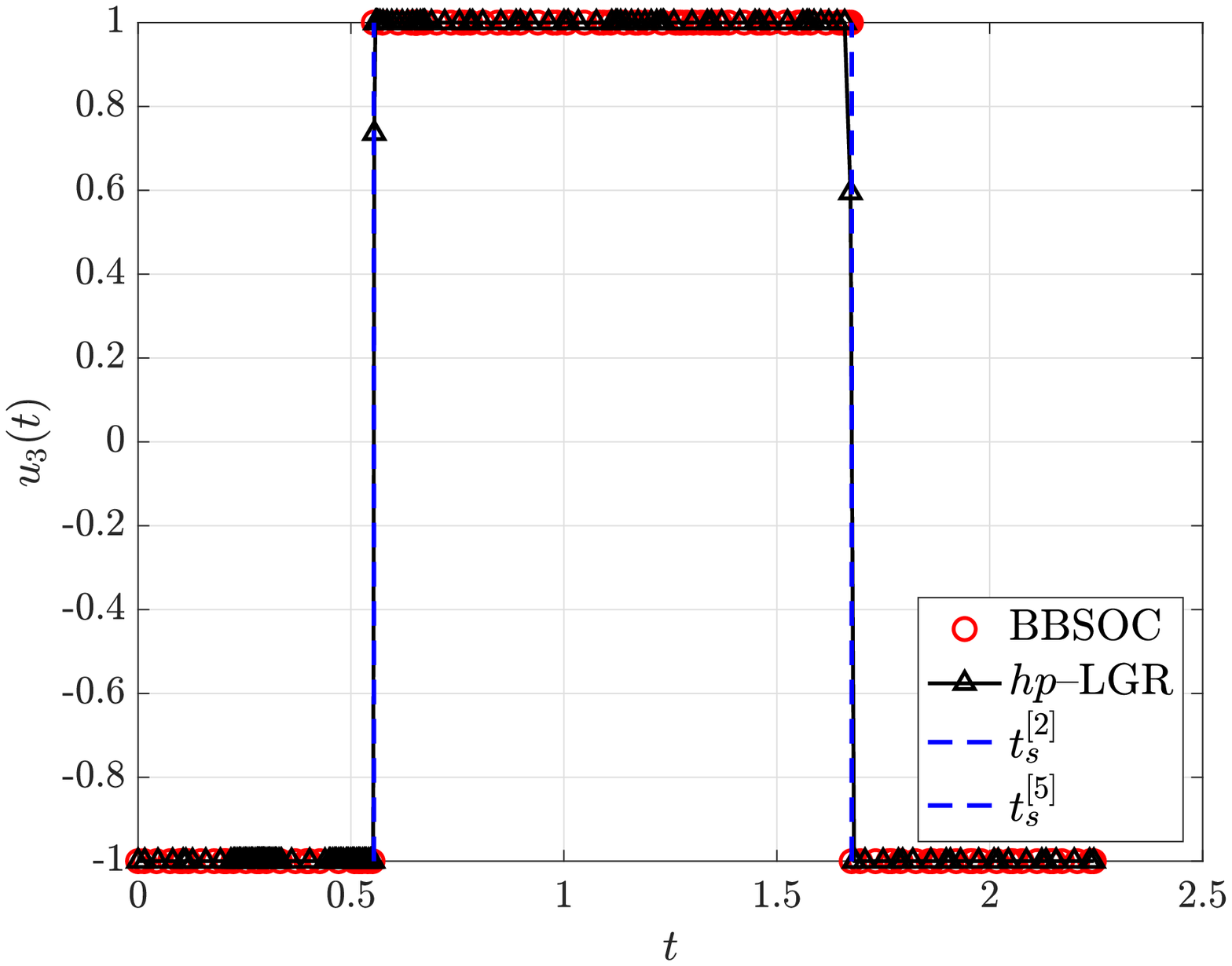}} \\ 
\caption{Control component solutions for the spinning NRTR maneuver. \label{fig:bangNRTR1}} 
\end{figure}

\begin{figure}[ht]
\centering
\vspace*{0.25cm}
\subfloat[\label{fig:bangNRTR-PHI}]{\includegraphics[scale=0.3]{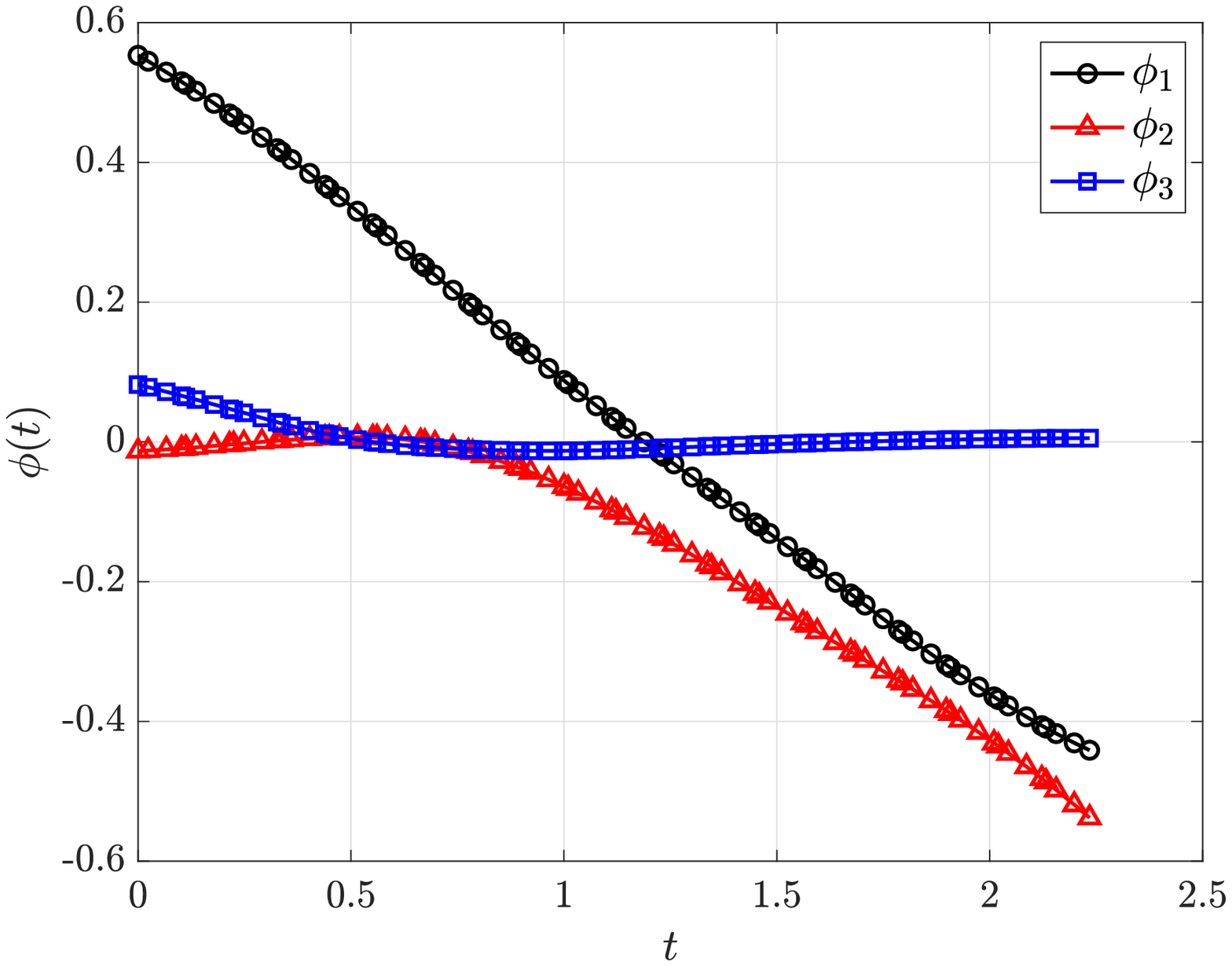}}~~\subfloat[\label{fig:bangNRTR-omega}]{\includegraphics[scale=0.3]{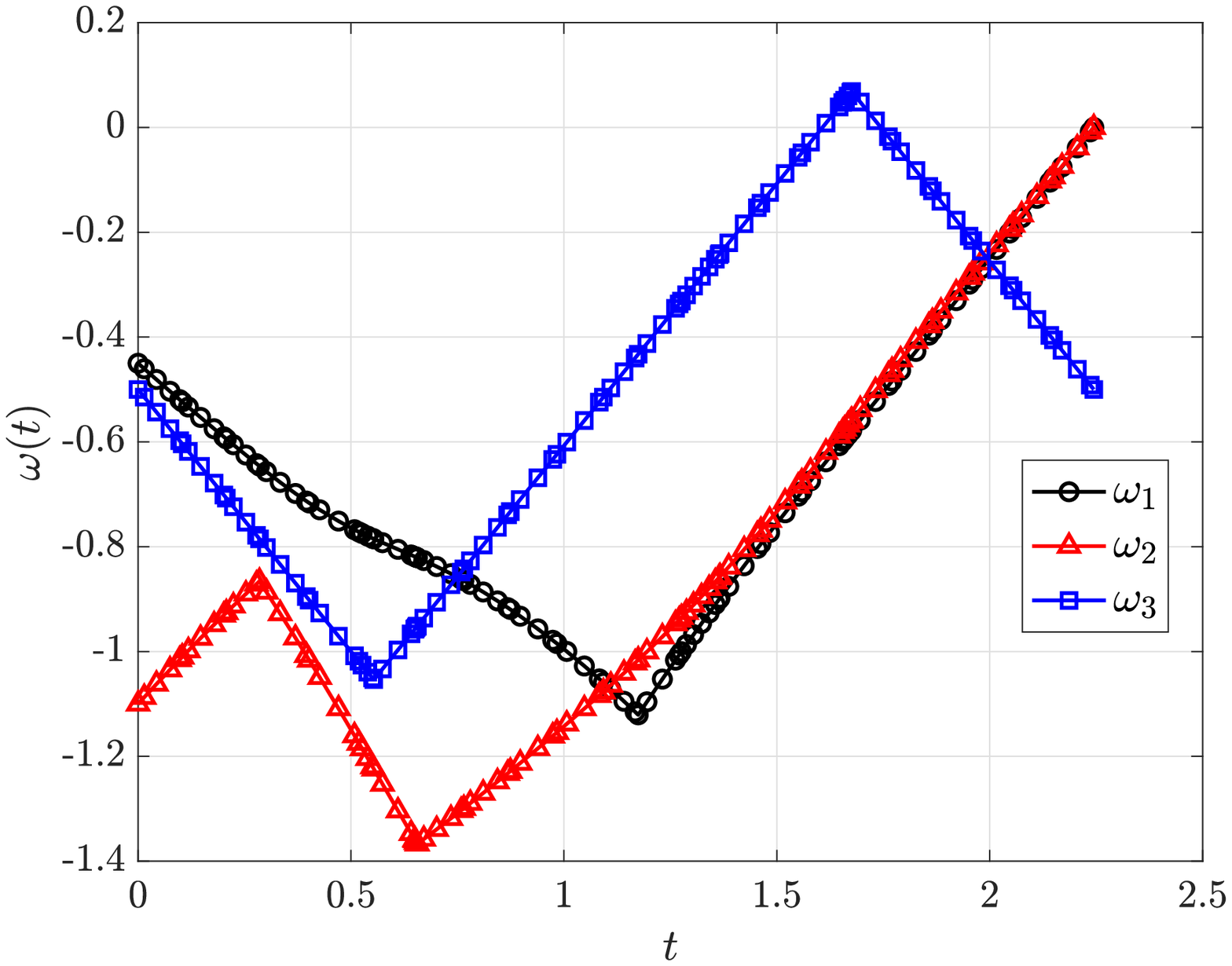}}\\
\caption{Switching functions and state solutions for the spinning NRTR maneuver. \label{fig:bangNRTR2}} 
\end{figure}

Figures~\ref{fig:bangNRTR1} and~\ref{fig:bangNRTR2} shows the three control solutions along with their corresponding switching functions and the time history of the angular velocities. Given that this scenario was not solved in Ref.~\cite{Shen1999} a direct comparison of results cannot be provided. However the two-torque control solution is obtained using the BBSOC method and produces a minimum time of $2.3069$ seconds as shown in Table~\ref{tab:bangNRTR}. Note that the switch  times cannot be compared as the switching structure is different between the two-torque and three-torque solutions. Using the three-torque control formulation results in a $2.70\%$ reduction in the minimum time. The values for the switch times and final time are found in Table~\ref{tab:bangNRTR} along with the CPU times for the BBSOC and $hp$-LGR methods. The BBSOC method converges faster than the $hp$--LGR method which is to be expected as the $hp$--LGR method often required more mesh refinements to converge.

\begin{table}[ht]
\caption{\label{tab:bangNRTR} Comparison of computational results for a NRTR maneuver.}
\centering
\begin{tabular}{ c " c c c c c c c c c c c } 
\hline \hline
                &  $t_{s}^{[1]}$  & $t_{s}^{[2]}$   & $t_{s}^{[3]}$ & $t_{s}^{[4]}$ & $t_{s}^{[5]}$ & $t_f$ & CPU~$[s]$   \\ \thickhline 
BBSOC  & $0.2851$  & $0.5536$ & $0.6557$  & $1.1745$ & $1.6759$  & $2.2445$  &  $2.24$   \\ 
$hp$--LGR  & $0.2882$  & $0.5548$  & $0.6630$  & $1.1742$ & $1.6731$  & $2.2445$  & $5.38$ \\
$2$--Torque & $-$ & $-$ & $-$ &$-$ & $-$ & $2.3069$ & $1.59$ \\
\hline \hline
\end{tabular}
\end{table}

\subsection{Non-Rest-to-Rest Maneuver for \texorpdfstring{$\omega_3=0$}{TEXT}: Finite Singular Arc\label{sect:finSingSoln-NRTR}}

For this special case the body is assumed to not be spinning about the symmetry axis, thus $(a,\omega_{3})=(0.5,0),\, \forall t\in[t_0,t_f]$. The boundary conditions
\begin{displaymath}
  \begin{array}{lclclcl}
    \omega_{10} & = & -0.45  &,&  \omega_{1f} & = & 0, \\
    \omega_{20} & = & -1.1  &,&  \omega_{2f} & = & 0,  \\
    \omega_{30} & = & 0 &,& \omega_{3f} & = & 0, \\
    x_{10} & = & 0.1 &,&  x_{1f} & = & 0,  \\
    x_{20} & = & -0.1 &,&  x_{2f} & = & 0,
  \end{array}
\end{displaymath}
are representative of a NRTR maneuver. An analysis of this special case was detailed in Section~\ref{sect:nonspinningCase} and showed that the resulting controls are bang-bang and singular. Specifically, the first control component contains a second-order singular arc and the second control component is bang-bang. Here the third control is not used as there is no rotation about the symmetry or third-axis. The optimal singular control is given by Eq.~\eqref{eq:FSA}. The optimal minimum time to complete the maneuver is $2.8839$ seconds.

\begin{figure}[ht]
\centering
\vspace*{0.25cm}
\subfloat[\label{fig:fsa-U1}]{\includegraphics[scale=0.3]{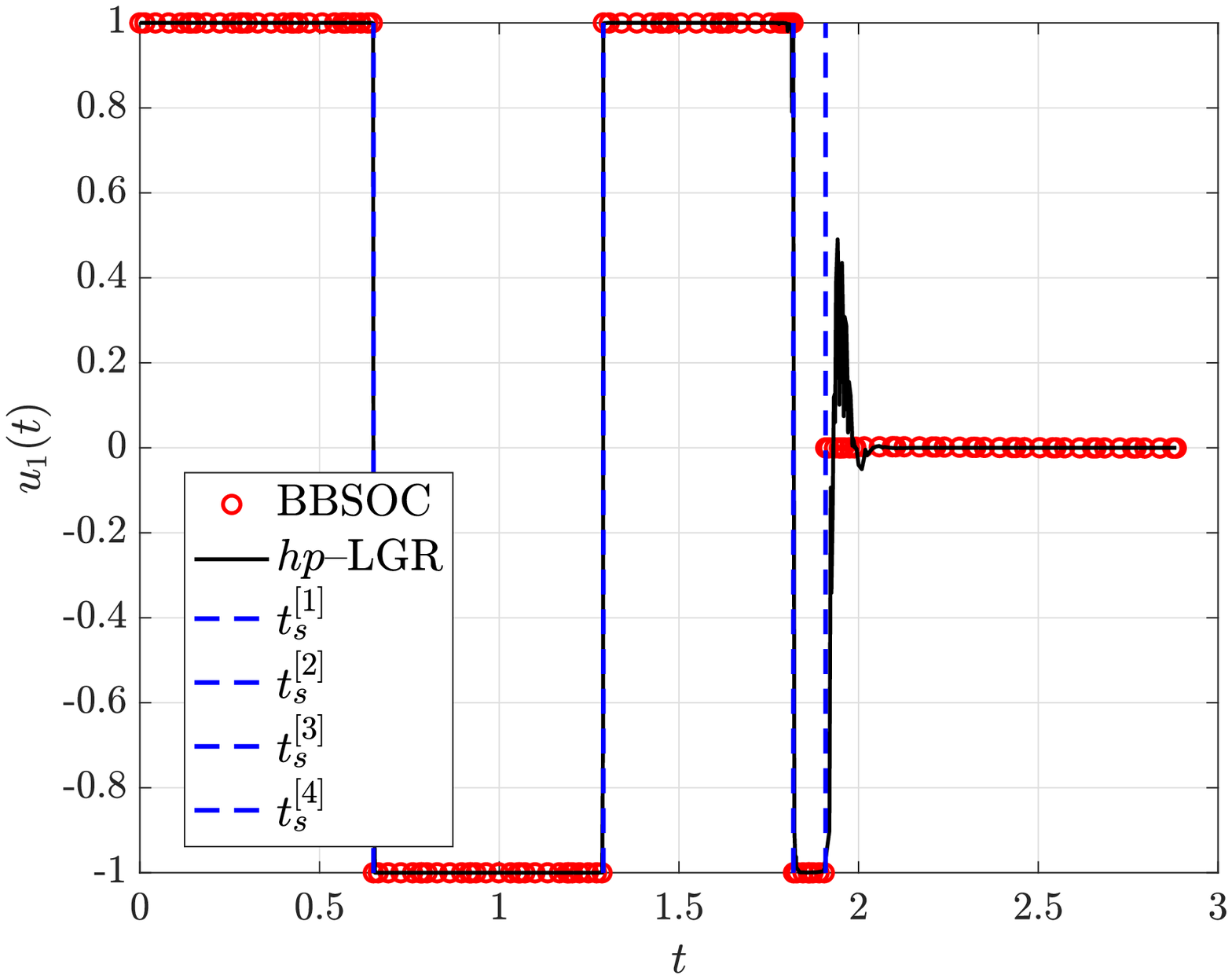}}~~\subfloat[\label{fig:fsa-U2}]{\includegraphics[scale=0.3]{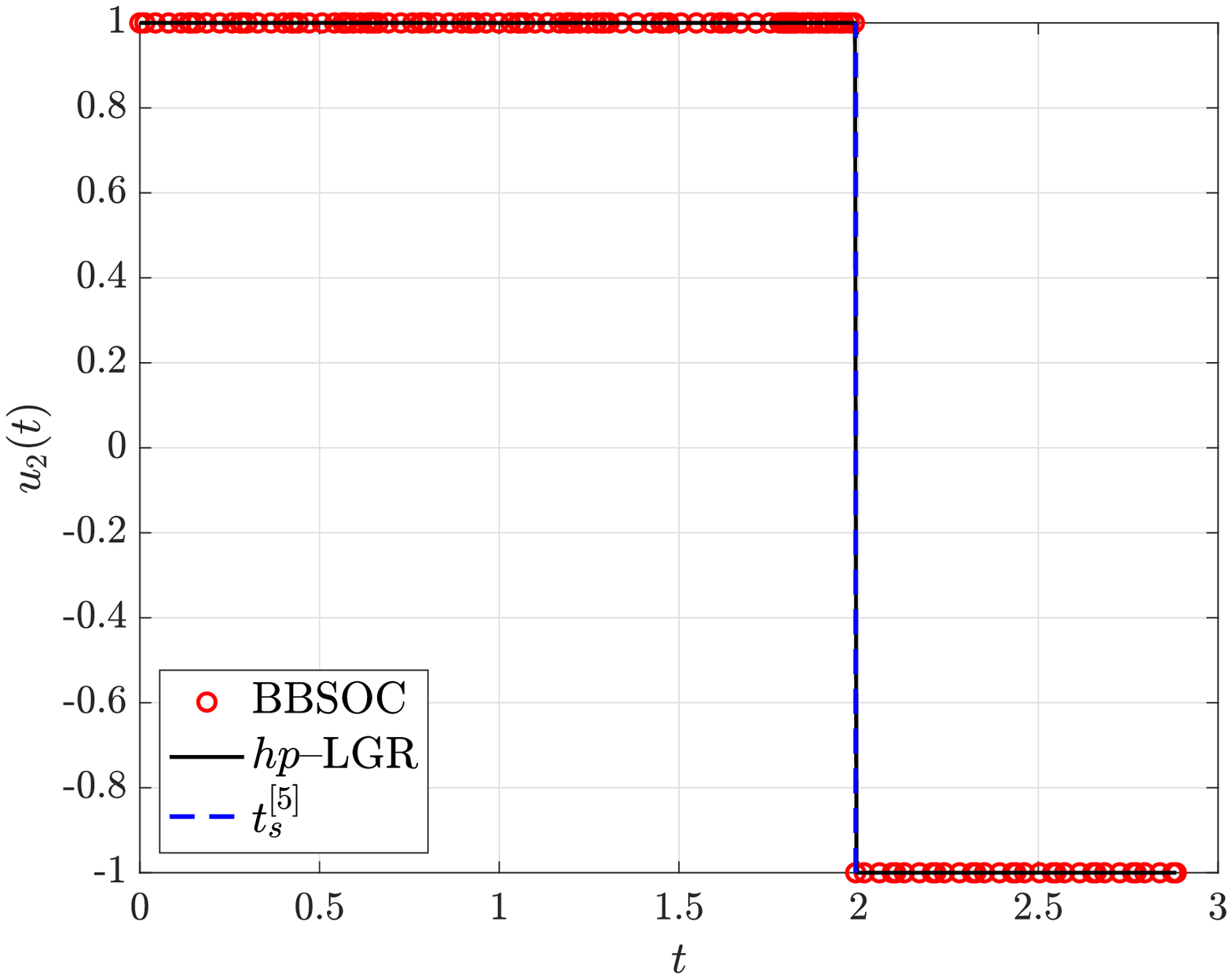}}~~\subfloat[\label{fig:fsa-Hist}]{\includegraphics[scale=0.3]{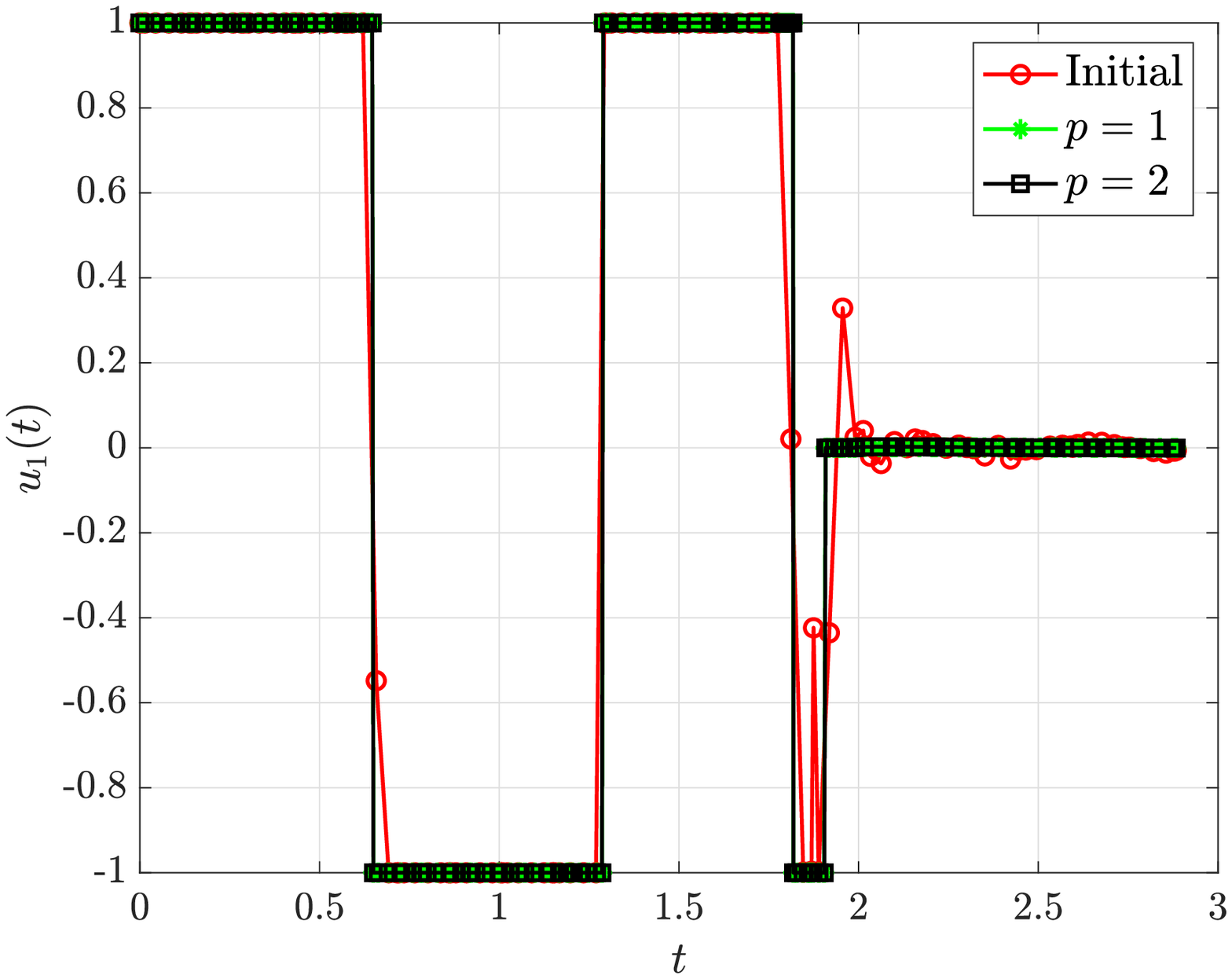}} \\
\caption{Control component solutions for the nonspinning NRTR maneuver, where Fig.~\ref{fig:fsa-Hist} shows the control history of the regularization procedure for obtaining the singular control.  \label{fig:fsa1}} 
\end{figure}

\begin{figure}[ht]
\centering
\vspace*{0.25cm}
\subfloat[\label{fig:fsa-PHI}]{\includegraphics[scale=0.3]{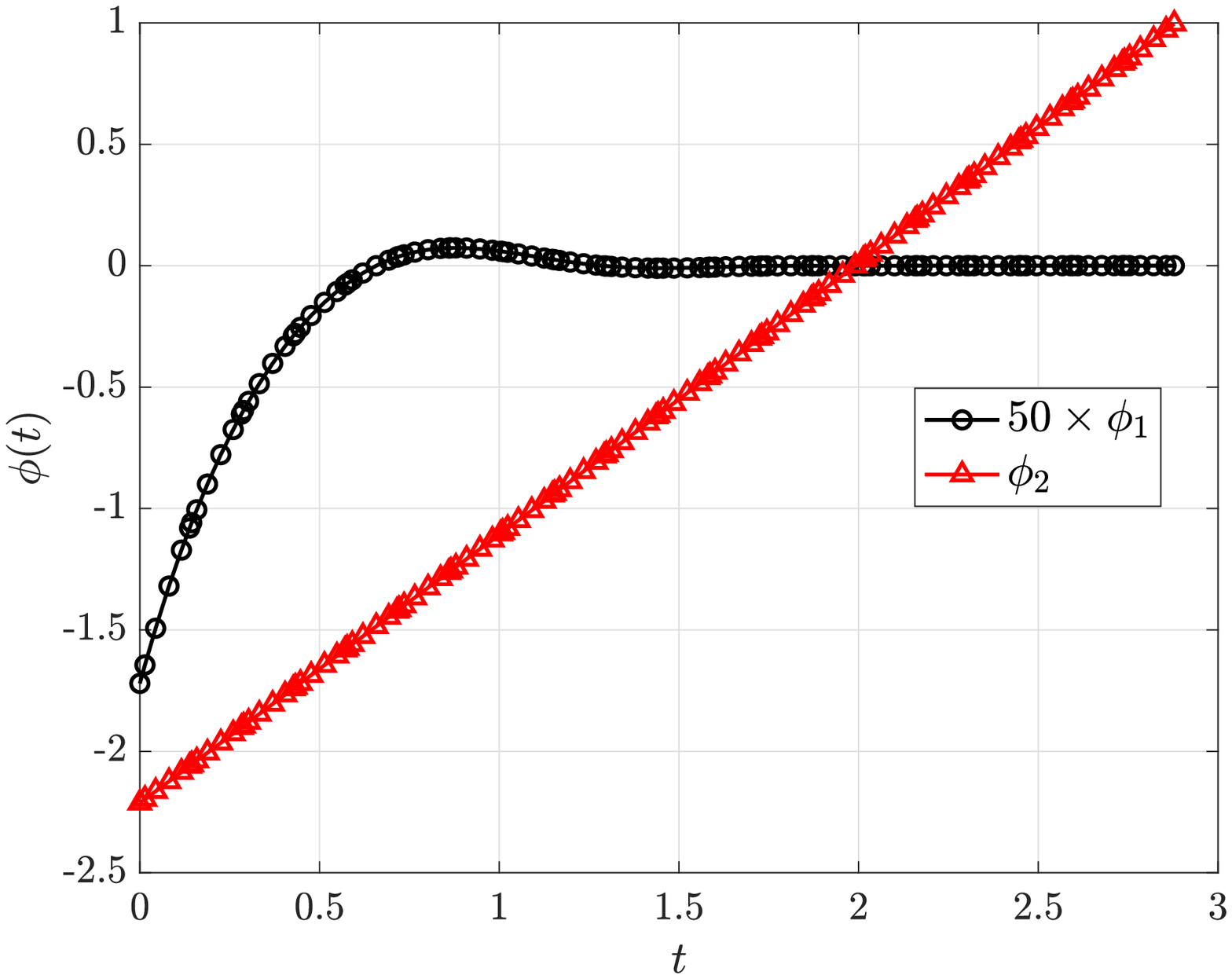}}~~\subfloat[\label{fig:fsa-omega}]{\includegraphics[scale=0.3]{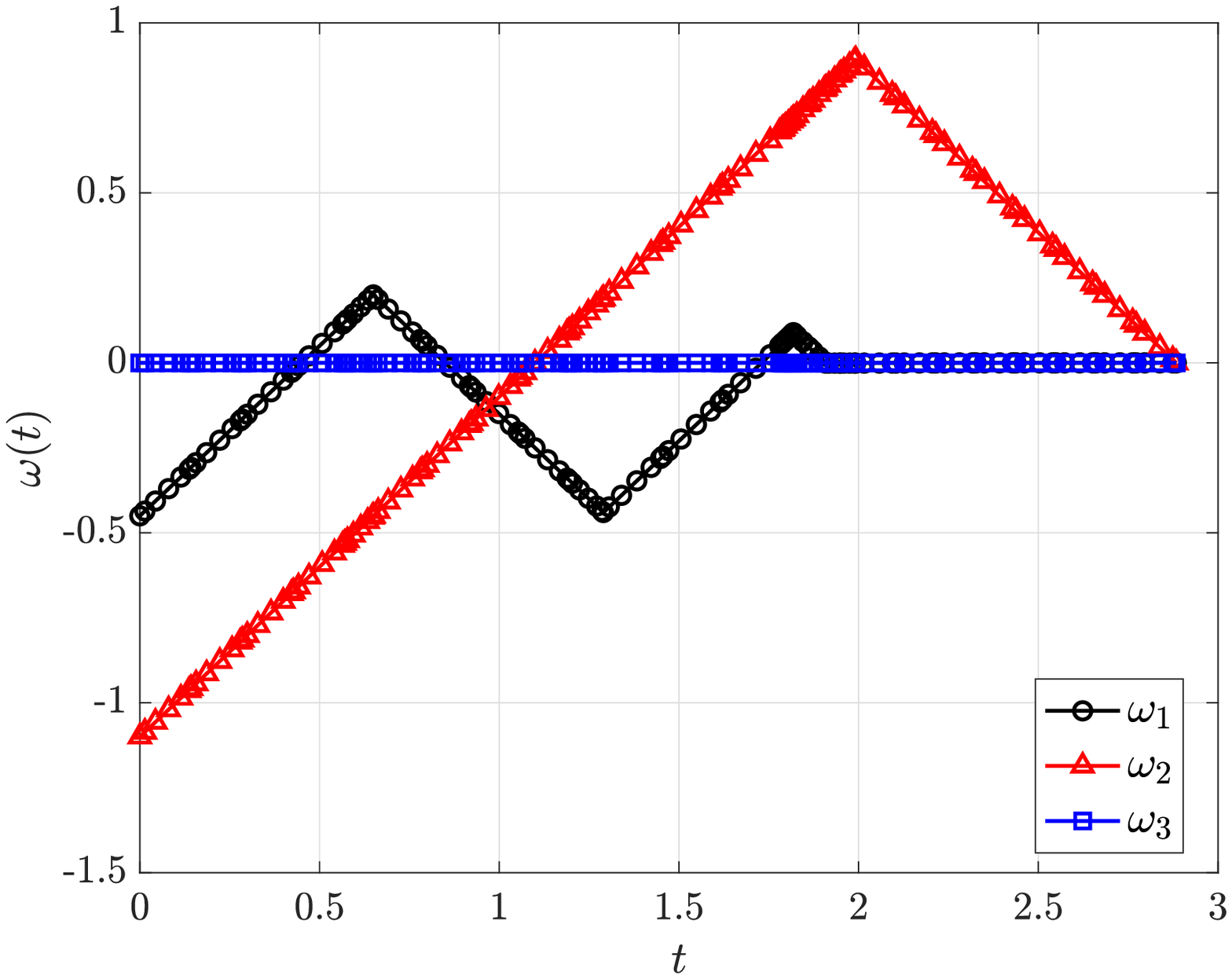}}~~\subfloat[\label{fig:fsa-H}]{\includegraphics[scale=0.3]{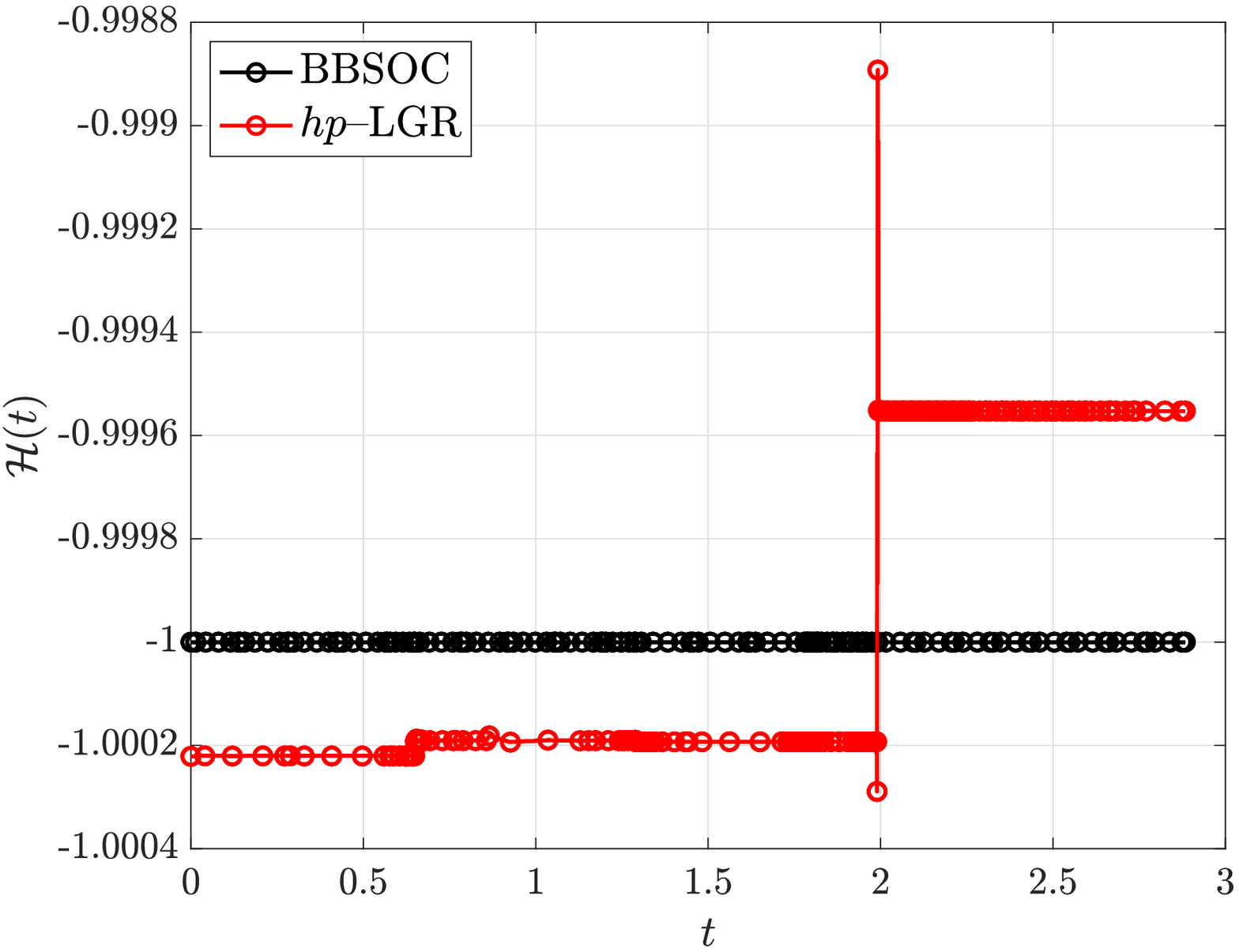}} \\ 
\caption{Switching functions and state solutions for the nonspinning NRTR maneuver, where Fig.~\ref{fig:fsa-H} shows the time history of the Hamiltonian. \label{fig:fsa2}} 
\end{figure}

Figures~\ref{fig:fsa1} and~\ref{fig:fsa2} show the numerical solution for the nonspinning NRTR maneuver. Figures~\ref{fig:fsa-U1}--\ref{fig:fsa-U2} show the optimal controls obtained using the BBSOC method and $hp$--LGR method. The singular arc begins at $t_s^{[4]}\approx 1.9054$. A numerical issue that occurs when using direct methods on singular optimal control problems is the occurrence of chattering behavior along the singular arc. This is observed in the $hp$--LGR solution and successfully removed by the BBSOC method as shown in Fig.~\ref{fig:fsa-Hist}. Figure~\ref{fig:fsa-Hist} shows the control history of the regularization procedure implemented by the BBSOC method. Additionally, values for the switch times and final time are provided in Table~\ref{tab:fsa} along with the CPU times for each method. Given that Ref.~\cite{Shen1999} does not provide the final time for their solution to more that two decimal places, it is difficult to compare the minimum times achieved by each of the numerical approaches. It is also noted that the CPU time listed for Ref.~\cite{Shen1999} is composed of the time it took for EZopt to converge to the initial guess for BNDSCO ($5$ minutes) and the computation time for BNDSCO ($2$ seconds). Regardless, it is clear that the BBSOC framework is more computationally efficient at solving nonsmooth and singular optimal control problems. 

Optimality is verified for the existence of the singular arc by checking that the generalized Legendre-Clebsch condition in Eq.~\eqref{eq:gen-Legendre-Clebsch} is satisfied and checking that the Hamiltonian satisfies the transversality condition in Eq.~\eqref{eq:transversalityCond}. According to Appendix~\ref{sect:appSpecial}~ $\omega_1=x_1=0$ for Eq.~\eqref{eq:gen-Legendre-Clebsch} to be satisfied. This is confirmed by looking at $\omega_1$ and $x_1$ in Figs.~\ref{fig:fsa-omega} and~\ref{fig:fsa-POS}. The second point can be verified by computing the Hamiltonian using the numerical solution obtained by each method. The Hamiltonian time history is provided in Fig.~\ref{fig:fsa-H} and shows that the BBSOC solution is optimal by producing a constant Hamiltonian value of $-1$ whereas the $hp$--LGR solution produces a nonconstant Hamiltonian.  The maximum absolute error of the Hamiltonian for both the BBSOC method and the $hp$--LGR method are $1.69\cdot 10^{-7}$ and $1.11\cdot 10^{-3}$, respectively. Thus, the singular arc is optimal based on the analysis of Section~\ref{sect:singArcAnalysis} and Appendix~\ref{sect:appSpecial}, and the optimality verification provided in this section.

\begin{table}[ht]
\caption{\label{tab:fsa} Comparison of computational results for the nonspinning NRTR maneuver.}
\centering
\begin{tabular}{ c " c c c c c c c c c c c c } 
\hline \hline
                &  $t_{s}^{[1]}$  & $t_{s}^{[2]}$   & $t_{s}^{[3]}$ & $t_{s}^{[4]}$ & $t_{s}^{[5]}$ & $t_f$ & $\delta$  & $\epsilon$  &  $p$ & CPU~$[s]$  \\ \thickhline 
BBSOC  & $0.6498$  & $1.2898$ & $1.8177$  & $1.9054$ & $1.9919$  & $2.8839$  &  $1.24 \times 10^{-10}$ & $10^{-3}$  &  $2$ & $3.37$ \\ 
$hp$-LGR  & $0.6498$  & $1.2897$  & $1.8184$  & $1.9226$ & $1.9922$  & $2.8839$  & $-$ & $-$ & $-$  & $5.18$ \\
Ref.~\cite{Shen1999}  & $-$  & $-$  & $-$  & $1.9040$ & $-$  & $2.8800$  & $-$ & $-$   & $-$  & $302$  \\
\hline \hline
\end{tabular}
\end{table}

\subsection{Rest-to-Non-Rest Maneuver for \texorpdfstring{$a=0$}{TEXT}: Infinite-Order Singular Arc\label{sect:infSingSoln-RTR}}
\begin{figure}[ht]
\centering
\vspace*{0.25cm}
\subfloat[\label{fig:isa-U1}]{\includegraphics[scale=0.3]{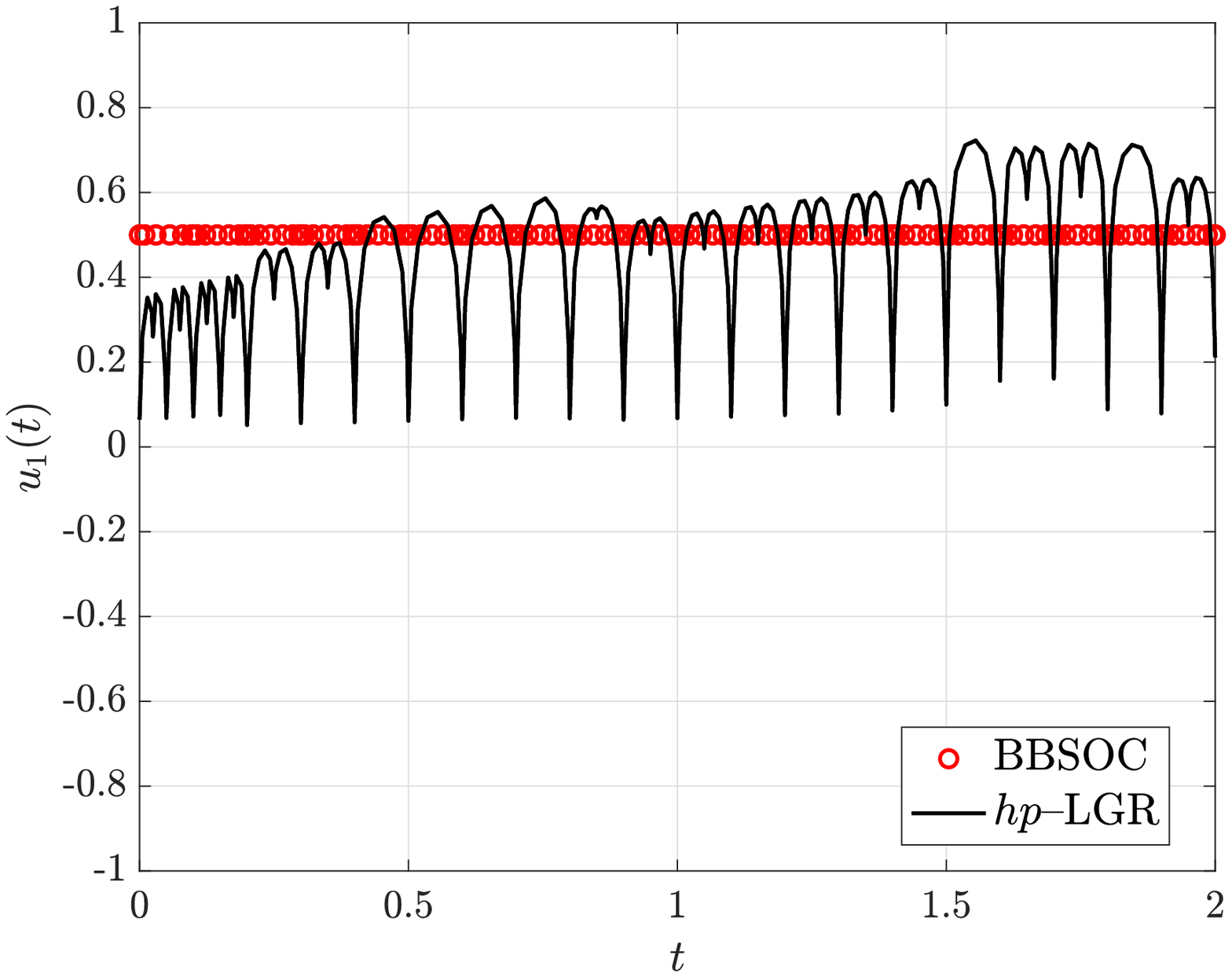}}~~\subfloat[\label{fig:isa-U2}]{\includegraphics[scale=0.3]{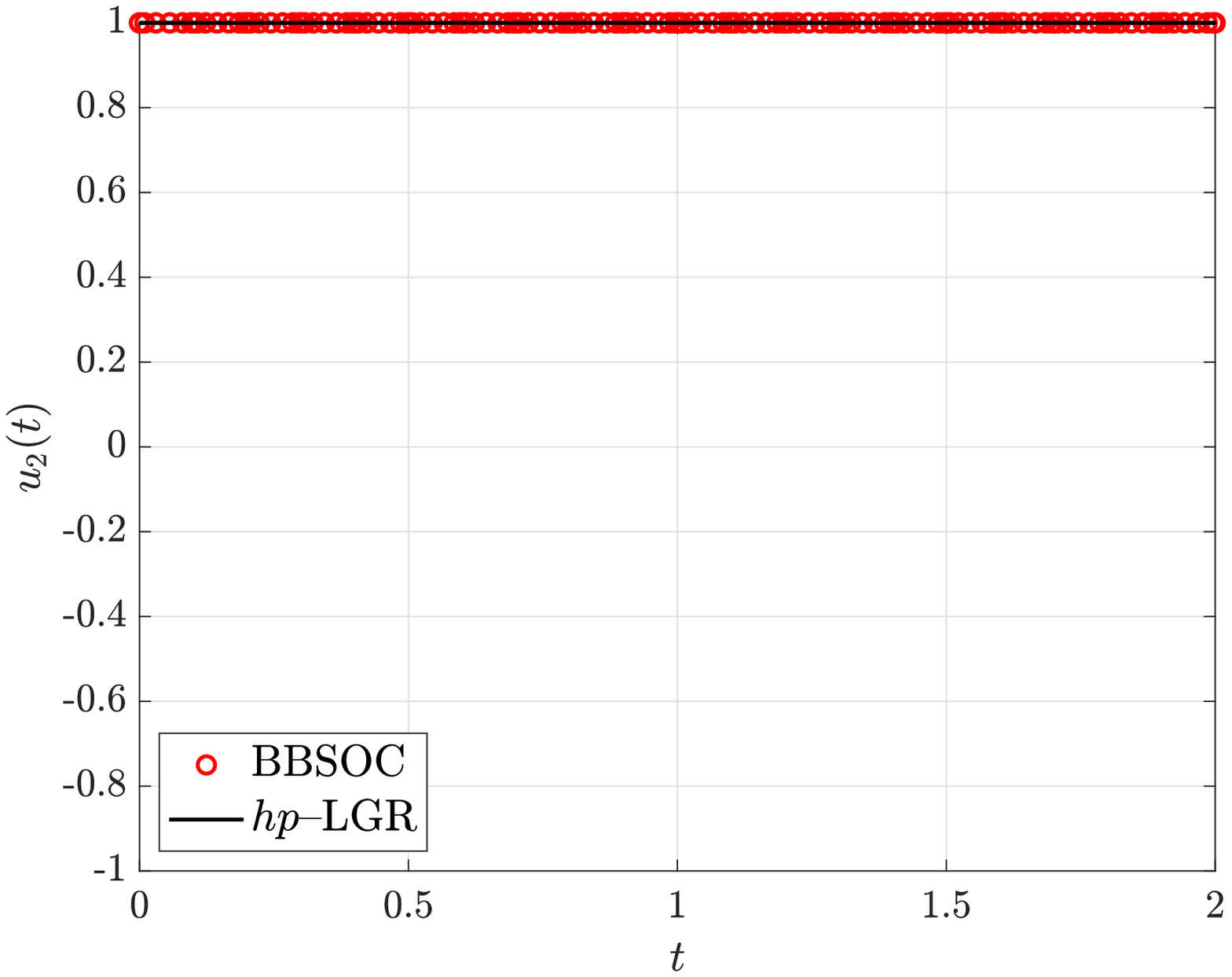}}~~\subfloat[\label{fig:isa-Hist}]{\includegraphics[scale=0.3]{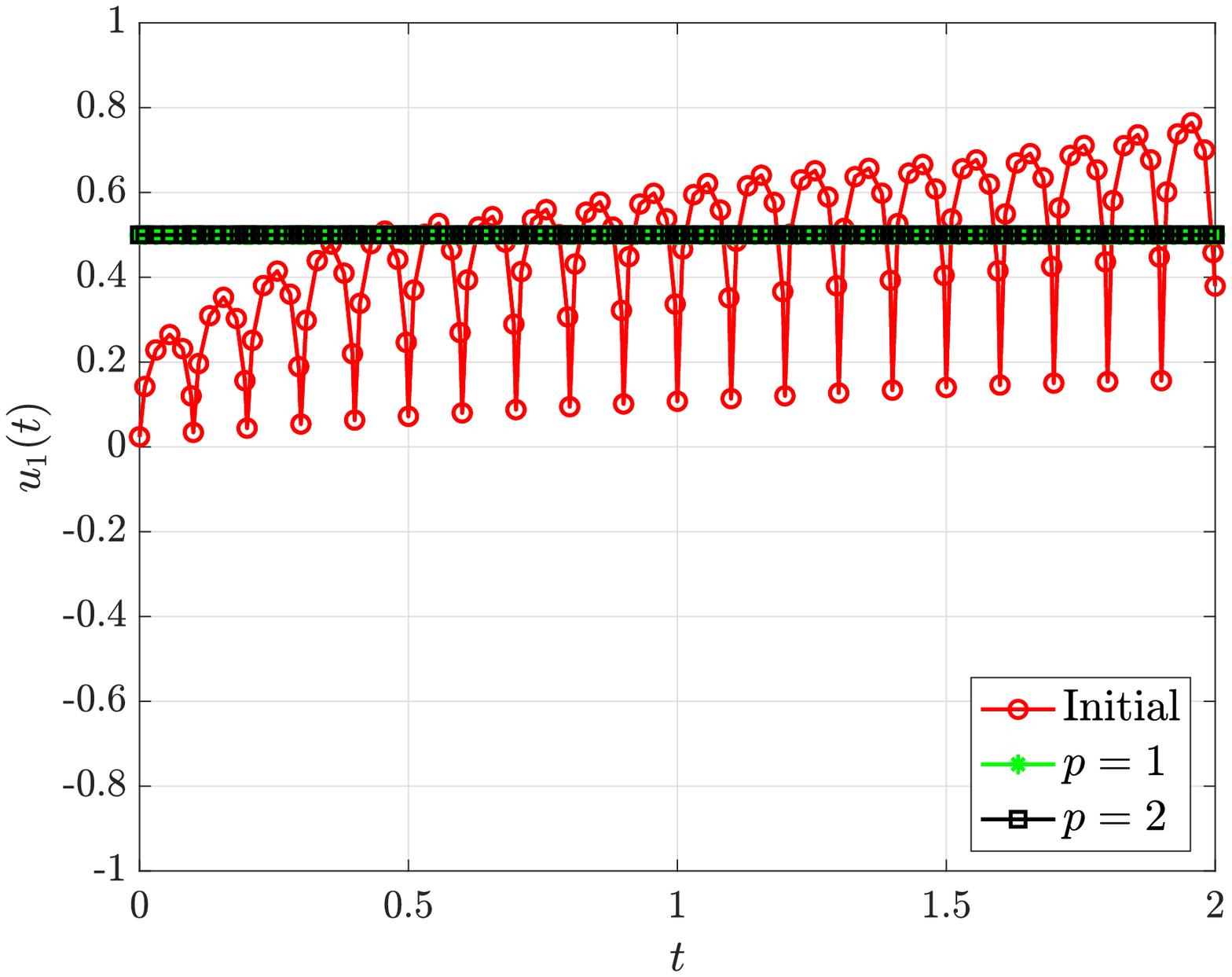}}
\caption{Control component solutions for the inertially symmetric RTNR maneuver, where Fig.~\ref{fig:fsa-Hist} shows the control history of the regularization procedure for obtaining the singular control.  \label{fig:isa1}} 
\end{figure}

\begin{figure}[ht]
\centering
\vspace*{0.25cm}
\subfloat[\label{fig:isa-PHI}]{\includegraphics[scale=0.3]{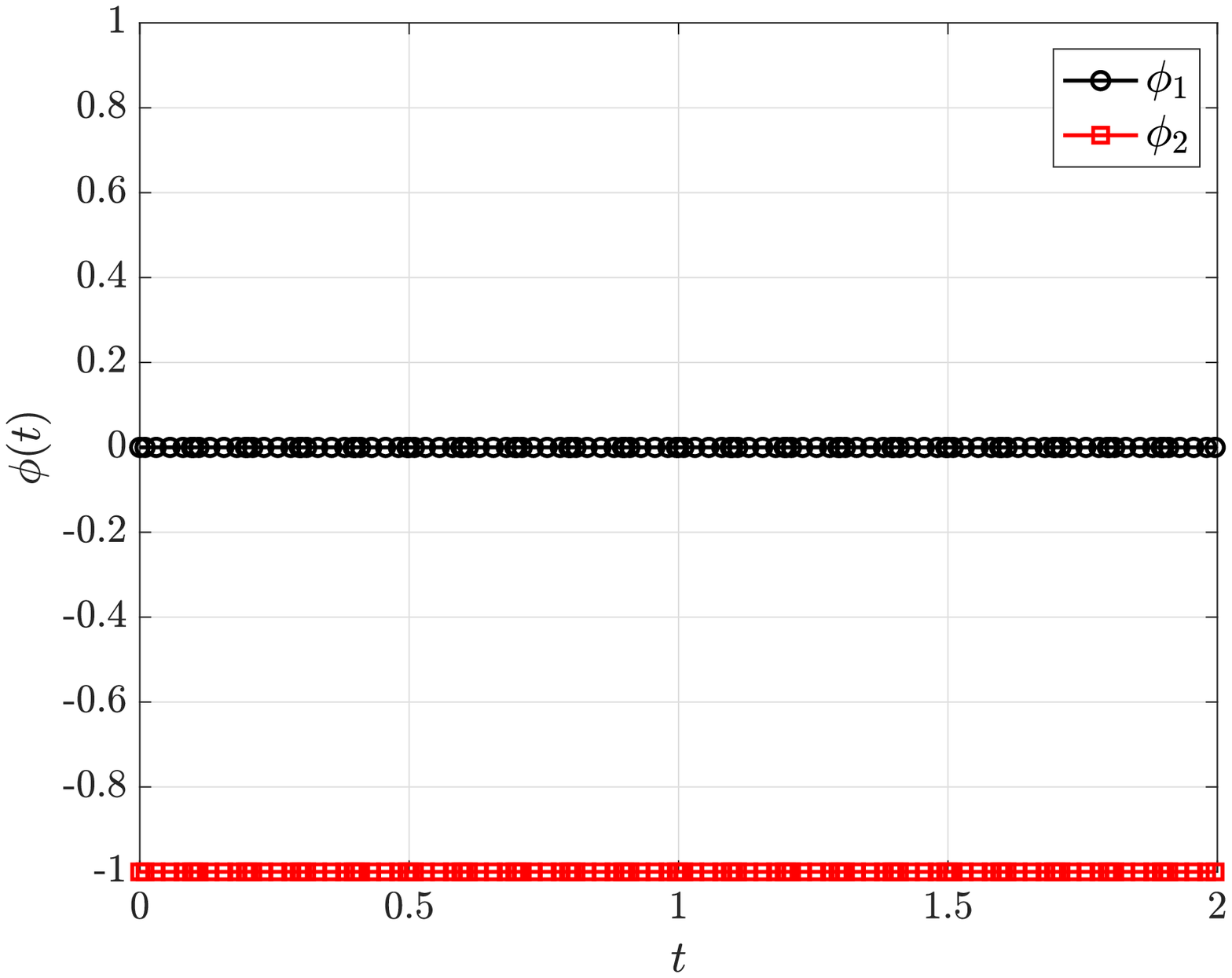}}~~\subfloat[\label{fig:isa-omega}]{\includegraphics[scale=0.3]{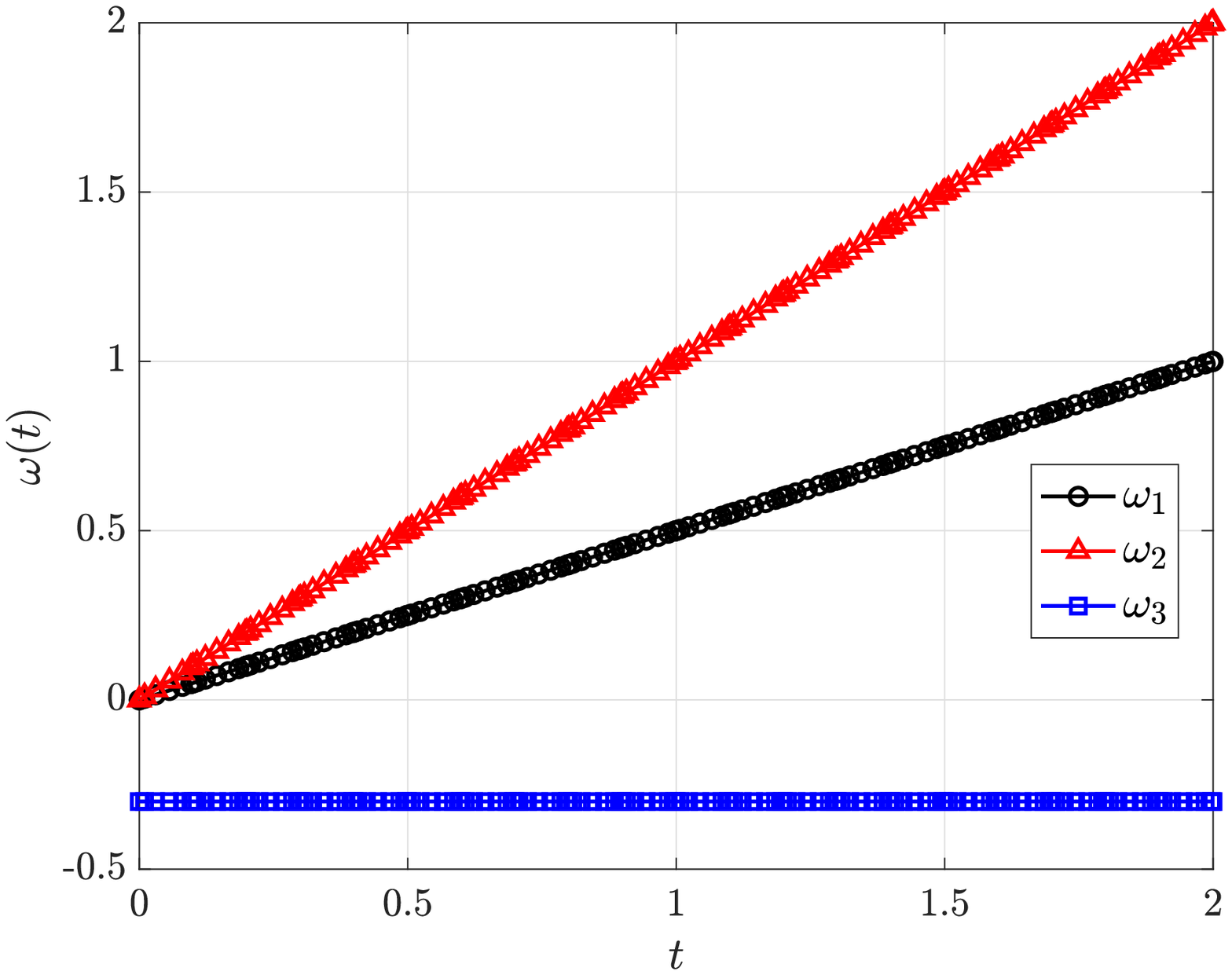}} \\
\caption{Switching functions and state solutions for the inertially symmetric RTNR maneuver. \label{fig:isa2}} 
\end{figure}

For this special case the rigid body is assumed to be inertially symmetric but spinning at a constant rate about the symmetry axis.  Consequently, $a=0$ and the boundary conditions are given as
\begin{displaymath}
  \begin{array}{lclclcl}
    \omega_{10} & = & 0 &,& \omega_{1f} & = & 1, \\
    \omega_{20} & = & 0 &,& \omega_{2f} & = & 2, \\
    \omega_{30} & = & -0.3 &,& \omega_{3f} & = & -0.3, \\
    x_{10} & = & 0 &,& x_{1f} & = & \textrm{Free},  \\
    x_{20} & = & 0 &,&  x_{2f} & = & \textrm{Free}, 
  \end{array}
\end{displaymath}
and represent a rest-to-non-rest (RTNR) maneuver.  The first control component contains an infinite-order singular arc across the entire time horizon, the second component is bang-bang, and the third control component is zero since $\omega_3$ is constant.  As a result of the infinite-order singular arc, many singular control solutions produce the optimal trajectory.  The fact that many singular solutions exist is explained by looking at the reduced dynamics in Eq.~\eqref{eq:IS-eom} where it is noted that $\omega_2$ must reach the specified final angular velocity of $2~\textrm{rad}/\textrm{s}$ in $2$ seconds, the optimal final time. One solution to the infinite-order singular arc is provided in Fig.~\ref{fig:isa-U1}. Likewise, the corresponding numerical solution is provided in Figs.~\ref{fig:isa-U2}--\ref{fig:isa-omega} including the second control component, the regularization iteration history for the singular control component, the corresponding control switching functions, and time histories of angular velocity. As expected, these results are in agreement with the results of Ref.~\cite{Shen1999}, see Table~\ref{tab:isa}, noting that the results of this work were obtained without using the known structure of control or the fact that the control lies on an infinite-order singular arc. The parameters associated with the BBSOC method are provided in Table~\ref{tab:isa} along with the CPU times for the BBSOC and $hp$--LGR methods.

Finally, optimality can be verified for the occurrence of the singular arc by checking that the  boundary conditions are satisfied and checking that the Hamiltonian satisfies the transversality condition in Eq.~\eqref{eq:transversalityCond}. The first point is verified by checking the state solution. The second point can again be verified by calculating the Hamiltonian using the numerical solution obtained by each method. The Hamiltonian time history is analyzed and shows that the maximum absolute error of the Hamiltonian for both the BBSOC method and the $hp$--LGR method are $9.99\cdot 10^{-9}$ and $1.61\cdot 10^{-8}$, respectively. Thus, the infinite-order singular arc is optimal based on the analysis of Section~\ref{sect:singArcAnalysis} and the optimality verification provided in this section.

\begin{table}[ht]
\caption{\label{tab:isa} Comparison of computational results for the spinning inertially symmetric RTNR maneuver.}
\centering
\begin{tabular}{ c " c c c c c c } 
\hline \hline
                      				&  $\C{J}^*$    & $\delta$  & $\epsilon$  & $p$ & CPU$[s]$  \\ \thickhline 
BBSOC       				   & $2.00$      &  $1.11 \times 10^{-19}$  & $10^{-1}$   &  $2$ & $1.77$ \\ 
$hp$-LGR                   & $2.00$   & $-$    & $-$  & $-$     & $1.74$  \\
Ref.~\cite{Shen1999}  & $2.00$     & $-$   & $-$   & $-$  & $-$  \\
\hline \hline
\end{tabular}
\end{table}

\section{Discussion}\label{sect:discussion}

The results of Section~\ref{sect:results} highlight several key aspects concerning the trajectory optimization of the rigid spacecraft reorientation problem described in Section~\ref{sect:probForm}. It was found that three-torque control formulation (as given in Ref.~\cite{Fleming2008}) produces an optimal trajectory with a smaller terminal time when compared with a two-torque control formulation because the spin rate about the symmetry axis is allowed to vary in the three-torque formulation. The lower terminal time obtained using the three-torque formulation is confirmed in the results of Sections~\ref{sect:bangSoln-RTR} and \ref{sect:bangSoln-NRTR}.  Comparisons cannot be made between the two-torque and three-torque control solutions when analyzing the special case maneuvers because the control torque about the symmetry axis is not used.

Next, Section~\ref{sect:results} explores the performance of the two numerical approaches discussed in Section~\ref{sect:methodComp}. Each of the maneuvers studied are solved using the BBSOC method, the $hp$--LGR method (using the implementation of Ref.~\cite{Patterson2014}), and the approach given in Ref.~\cite{Shen1999} (excluding the NRTR maneuver in Section~\ref{sect:bangSoln-NRTR}). For the bang-bang maneuvers, the BBSOC method outperforms the method of Ref.~\cite{Shen1999} both in computation time and accuracy. For the special case maneuvers where finite and infinite-order singular arcs occur in the solution, again the BBSOC method outperforms the other methods mentioned. The key observation here is that the BBSOC method obtains more accurate solutions to these nonsmooth and singular problems and does not require that the singular arc conditions be enforced whereas the singular arc conditions must be enforced when using indirect shooting.  Furthermore, the BBSOC method neither requires accurate initial guesses of the solution nor requires a guess of the switching structure.  In contrast, the results of the method in Ref.~\cite{Shen1999} are only obtained if BNDSCO is provided a near optimal initial guess and the singular arc conditions are  provided.

Optimality is also verified for the maneuvers studied in Section~\ref{sect:results}. For the special case maneuvers in Sections~\ref{sect:finSingSoln-NRTR}--\ref{sect:infSingSoln-RTR} it was shown that the Hamiltonian is constant and numerically equal to $-1$.  In contrast, the solutions obtained using the $hp$--LGR method had numerical error on the order of $10^{-3}$ and $10^{-8}$, respectively. The generalized Legendre-Clebsch condition~\eqref{eq:gen-Legendre-Clebsch} was also shown to be satisfied for both singular controls. The optimality of the bang-bang solutions are also verified by considering the maximum absolute numerical error accrued in the Hamiltonian. For the bang-bang RTR maneuver the error in the Hamiltonian was found to be $1.50\cdot10^{-6}$ and $1.98\cdot10^{-3}$ for the BBSOC and $hp$--LGR methods, respectively.  For the bang-bang NRTR maneuver the error in the Hamiltonian was found to be $1.74\cdot10^{-6}$ and $2.58\cdot10^{-3}$ for the BBSOC and $hp$--LGR methods, respectively. 

\begin{figure}[ht]
\centering
\vspace*{0.25cm}
\subfloat[\label{fig:bangRTR-POS}]{\includegraphics[scale=0.3]{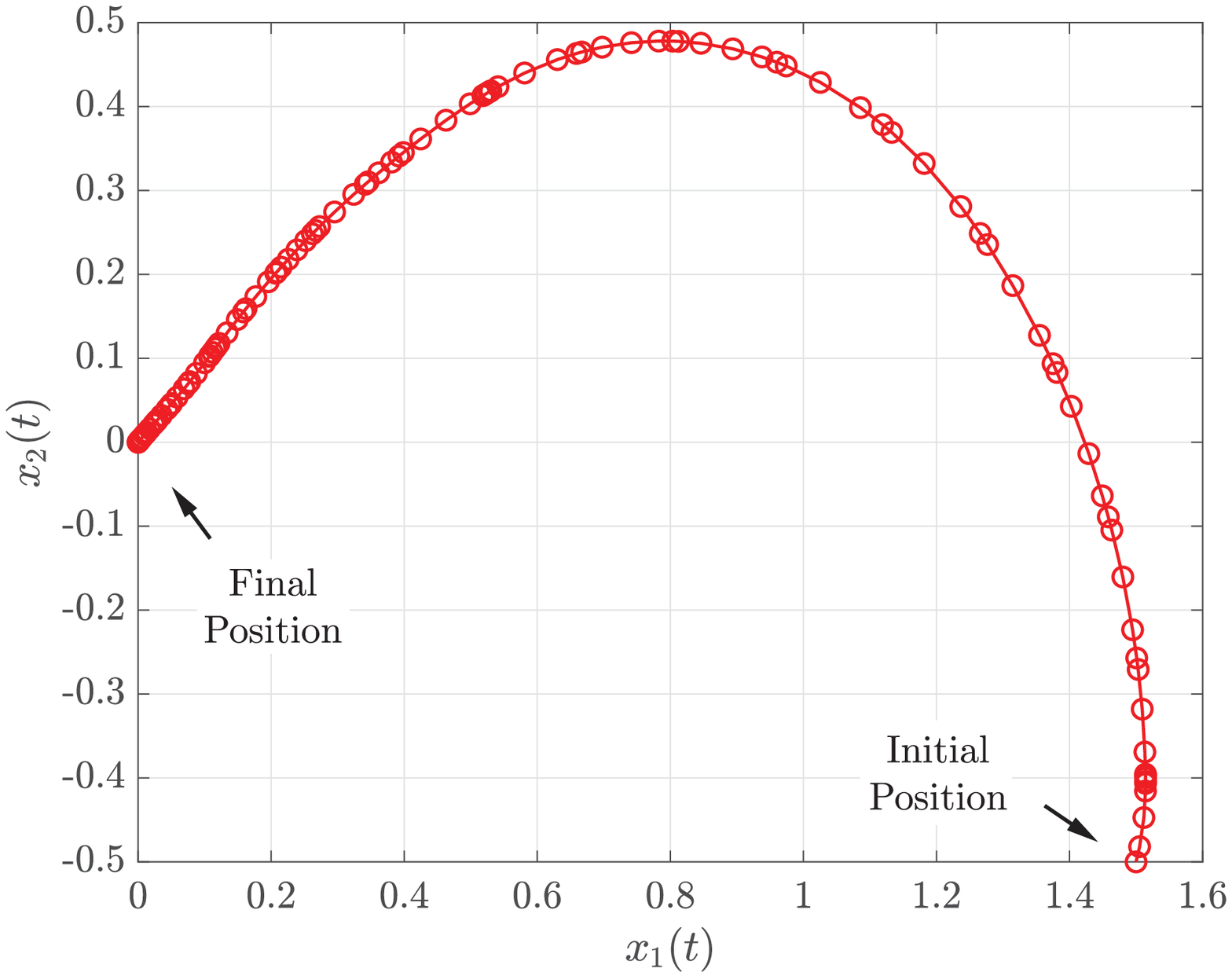}}~~\subfloat[\label{fig:bangNRTR-POS}]{\includegraphics[scale=0.3]{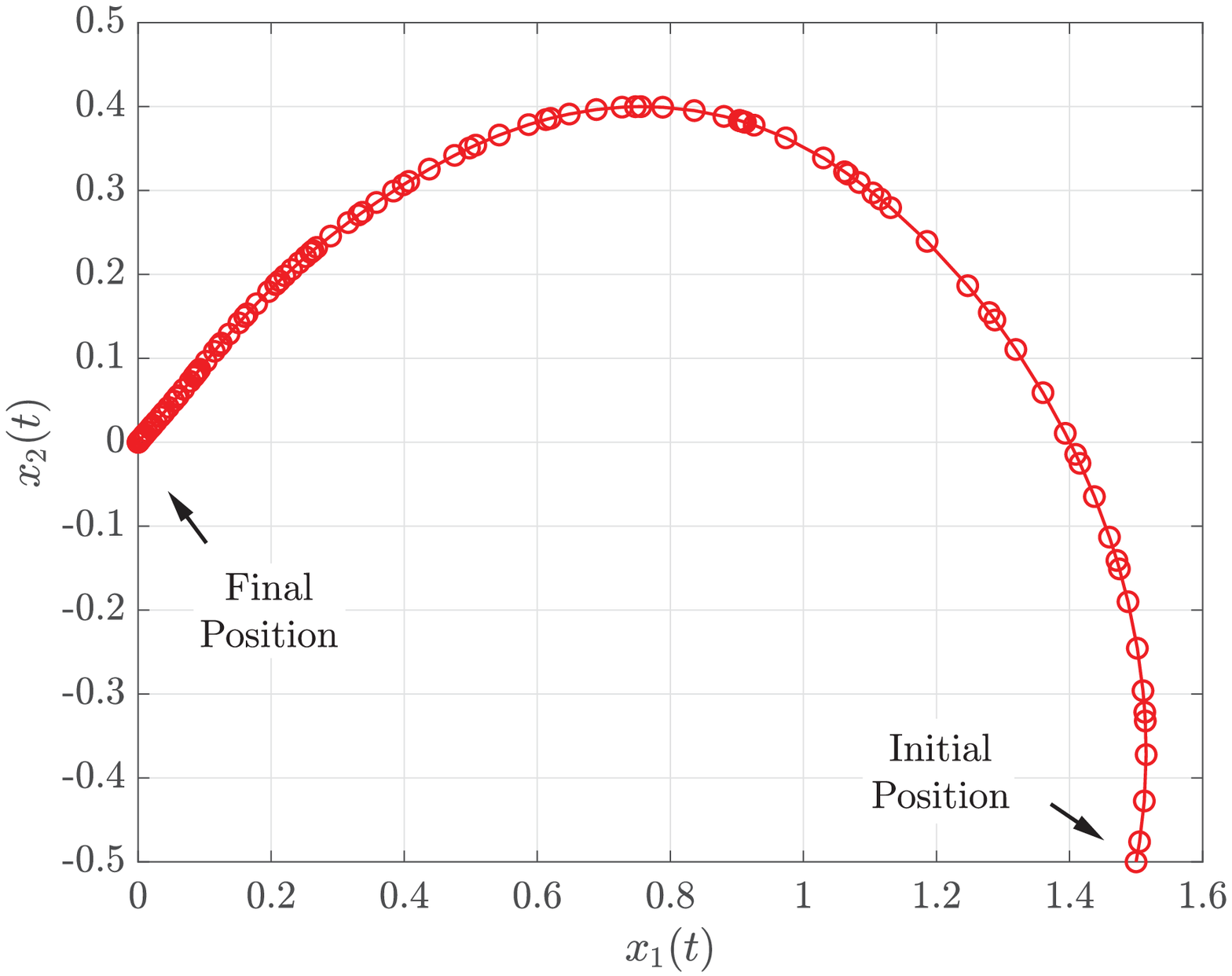}} \\
\subfloat[\label{fig:fsa-POS}]{\includegraphics[scale=0.3]{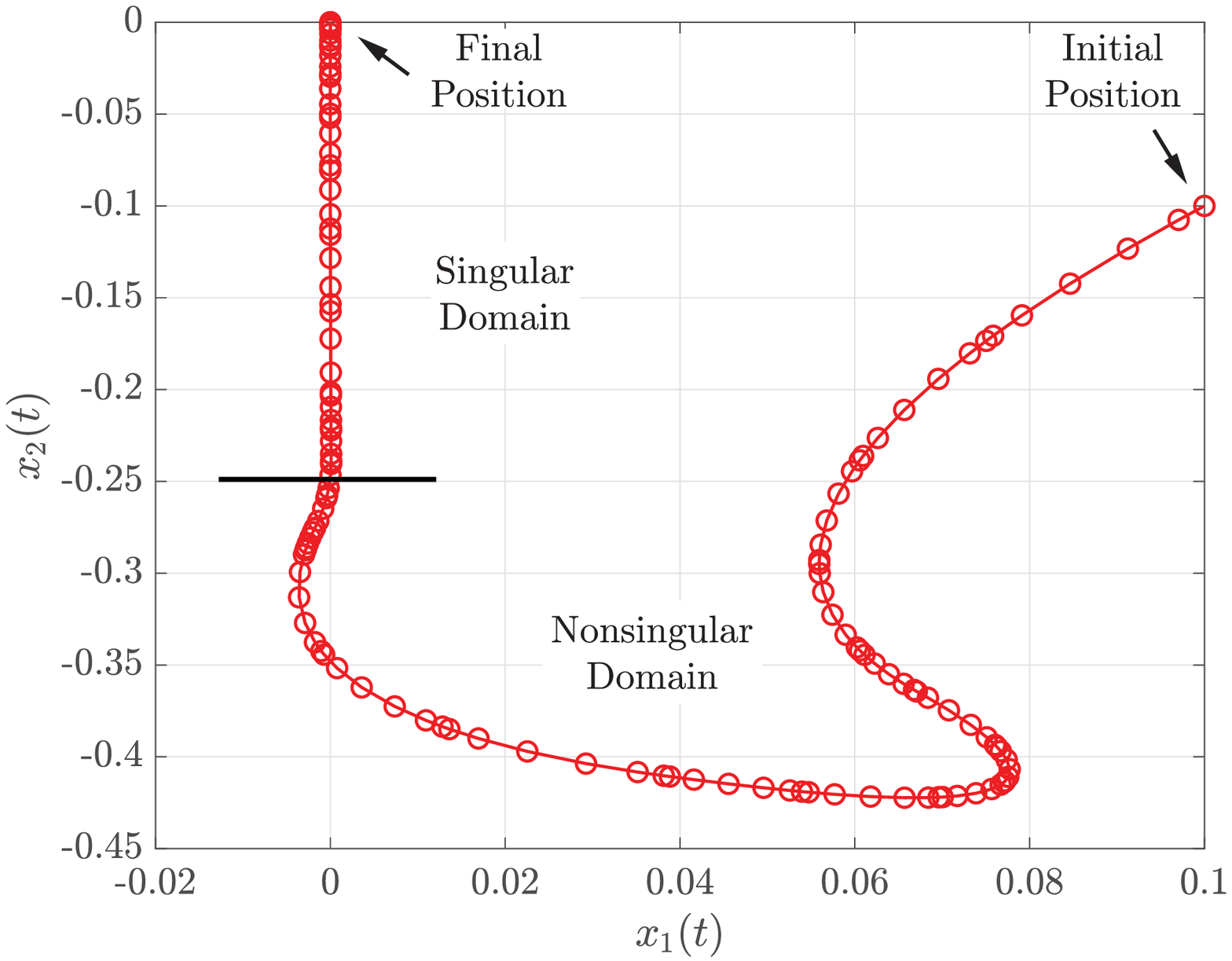}}~~\subfloat[\label{fig:isa-POS}]{\includegraphics[scale=0.3]{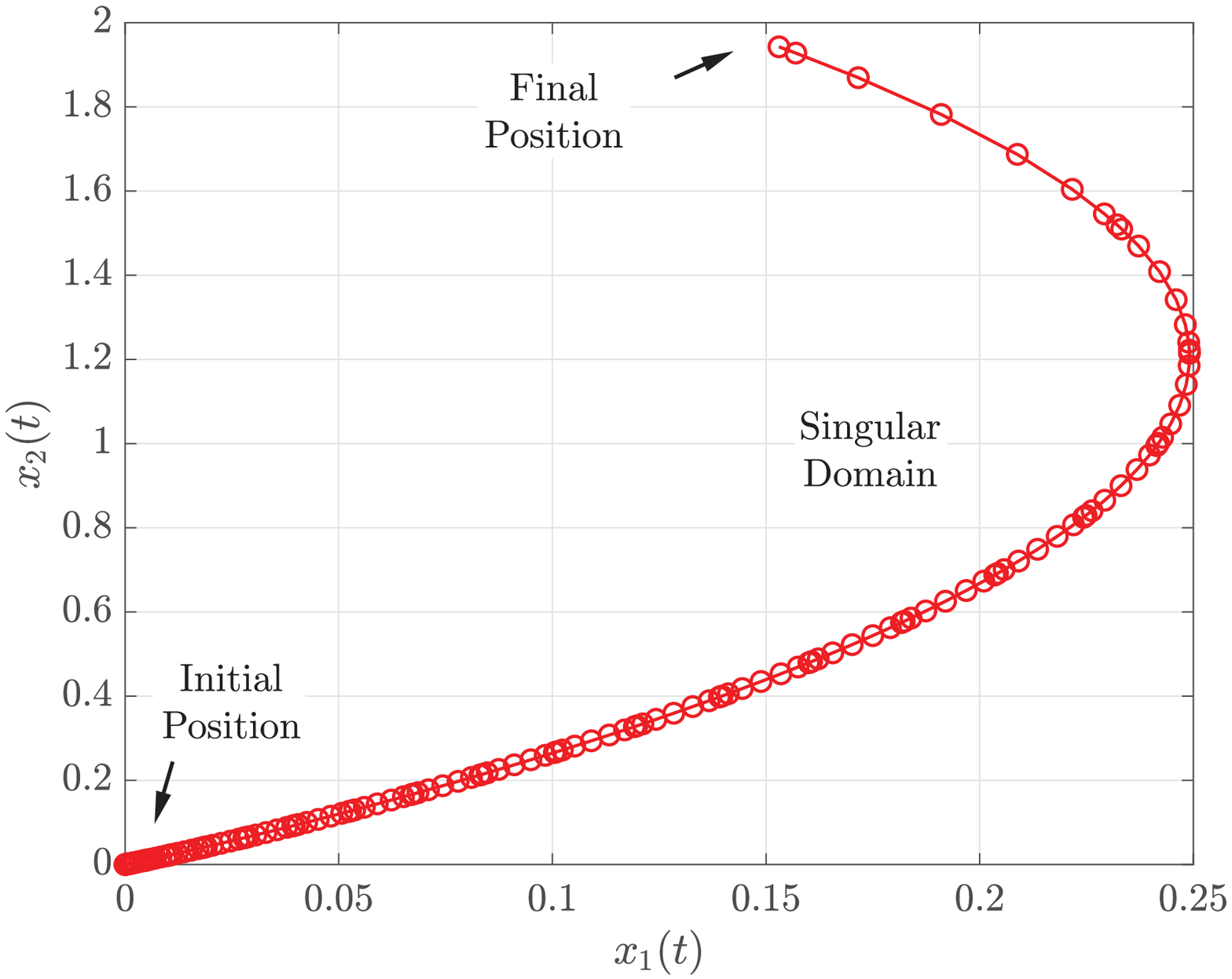}} \\
\caption{Optimal trajectories for the (\ref{fig:bangRTR-POS}) RTR (\ref{fig:bangNRTR-POS}) NRTR (\ref{fig:fsa-POS}) nonspinning NRTR and (\ref{fig:isa-POS}) inertially symmetric RTNR maneuvers in the $x_1-x_2$ plane.} 
\label{fig:trajectory}
\end{figure}

Finally, the optimal trajectories for each maneuver are discussed. Figure~\ref{fig:trajectory} shows the optimal trajectories of each maneuver in the $x_1-x_2$ plane representing the goal of reorienting the relative position of the inertial axis $\m{n}_3$, with respect to the body-fixed frame $\m{b}$, from an initial position to a final position. In particular, Fig.~\ref{fig:fsa-POS} shows the trajectory for the NRTR maneuver when $\omega_3=0$ in Section~\ref{sect:finSingSoln-NRTR}. Here the straight line portion of the trajectory represents an eigenaxis rotation when the rigid spacecraft is not spinning about the symmetry axis. An eigenaxis rotation does not appear in any of the other solutions since eigenaxis rotations in general are not time-optimal as shown in Refs.~\cite{Shen1999,Seywald1993}, but they do often produce the shortest angular trajectory between the rigid spacecrafts current attitude and the desired attitude. For this special case, however, the singular arc control is equivalent to an eigenaxis rotation along that portion of the trajectory and in this case is optimal. 

\section{Conclusions\label{sect:conc}}
The minimum-time reorientation of an axisymmetric rigid spacecraft under the influence of three control torques has been considered. The problem first considered by Ref.~\cite{Shen1999} is revisited to show improvement of numerical results using new methods. A review and further insight on the mathematical analysis of the necessary conditions for optimality has been carried out. Different maneuvers, including special cases, of the rigid spacecraft reorientation problem have been solved and analyzed using a recently-developed direct collocation method, called the BBSOC method, that is designed to solve bang-bang and singular optimal control problems..  Comparisons to a previous study of the problem of interest were made and show the key aspects of using the BBSOC method over previously developed methods. Overall the performance of the BBSOC method has proven to be more efficient and successful when compared to previous results in the literature. Although necessary conditions for optimality of the singular controls were developed and provided, they were not needed by the BBSOC method to obtain accurate solutions. In contrast, the numerical approach of Ref.~\cite{Shen1999} required the enforcement of these necessary conditions in their numerical approach. Additionally, the results have shown that considering three-torque control over two-torque control leads to optimal minimum-time trajectories with smaller terminal times.


\section{Acknowledgments}\label{sect:ack}

The authors gratefully acknowledge support for this research from the U.S.~National Science Foundation under grants DMS-1819002 and CMMI-2031213, the U.S.~Office of Naval Research under grant N00014-19-1-2543, and from the U.S.~Air Force Research Laboratory under contract FA8651-21-F-1041.  

\appendix

\section{Appendix}\label{sect:appendix}
\subsection{Derivation of Singular Control for \texorpdfstring{$a\neq0$}{TEXT} and \texorpdfstring{$\omega_3 \neq 0$}{TEXT}}\label{sect:appGeneral}

Consider the case where $u_1$ is singular and both $u_2$ and $u_3$ are bang-bang.  On an interval where $u_1$ is singular, the switching function $g_1$ must be differentiated four times with respect to time in order to determine the singular control.  First, it is known from Eq.~\eqref{eq:ham} that
\begin{equation}\label{g1=0:u1-singular}
  g_1 = \lambda_1 = 0
\end{equation}
when $u_1$ is singular.  Differentiating $g_1$ with respect to time four times yields the following results:
\begin{eqnarray}\label{eq:dg1dt}
  \frac{d g_1}{dt} & = & \lambda_{2}a\omega_{3} - \lambda_{4}\frac{1+x_1^2-x_2^2}{2} -\lambda_{5}x_1x_2, \\\label{eq:d2g1dt2}
  \frac{d^2g_1}{dt^2} &  = & \omega_{3} (1+a)\dot{\lambda}_2 + \omega_2(\lambda_3 x_2 - \lambda_4 x_1), \\\label{eq:d3g1dt3}
  \frac{d^3g_1}{dt^3} & = & a\omega_{3}^3(1+a)\lambda_2 - \omega_{3}\omega_1(1+2a)(\lambda_3x_2-\lambda_4x_1) + \dot{\lambda}_2\omega_1\omega_2 + a\omega_{3}\omega_2^2\lambda_2 + u_2(\lambda_3x_2-\lambda_4x_1), \\
  \frac{d^4g_1}{dt^4} & = & \nonumber \dot{\lambda}_2[ 2a\omega_{3}^3(a+1)^2 - 2a\omega_{3}\omega_1^2 + 2a\omega_{3} \omega_2^2 + 2\omega_1u_2 ] - \omega_2\lambda_2(4a^2\omega_{3}^2\omega_1 - 3a\omega_{3}u_2 ) \\ \label{eq:d4g1dt4}
  {} & + & \left[\dot{\lambda}_2 \omega_2 - \omega_{3}(1+2a)(\lambda_3x_2 - \lambda_4x_1) \right]  u_1, 
\end{eqnarray}  
where Eq.~\eqref{eq:lam-1} has been used in Eq.~\eqref{eq:dg1dt}, the fact that $\lambda_1=0$ from Eq.~\eqref{g1=0:u1-singular} has been used in Eq.~\eqref{eq:d2g1dt2}, the facts that $\lambda_1 = 0$ and $\dot{\lambda}_1 = 0$ from Eqs.~\eqref{g1=0:u1-singular} and \eqref{eq:dg1dt}, respectively, have been used in Eq.~\eqref{eq:d3g1dt3}, and the expression for $\dot{\lambda}_2$ from Eq.~\eqref{eq:lam-2} has been used in Eq.~\eqref{eq:d4g1dt4}.  
Note that all four derivatives in Eqs.~\eqref{eq:dg1dt}--\eqref{eq:d4g1dt4} are zero on a singular arc because $g_1$ is itself zero on the singular arc.  Therefore, $d^4g_1/dt^4$ in Eq.~\eqref{eq:d4g1dt4} is zero on a singular arc and the singular control $u_1$ is obtained as
\begin{equation}\label{eq:singularU1}
  u_1 = -\frac{\dot{\lambda}_2[ 2a\omega_{3}^3(a+1)^2 - 2a\omega_{3}\omega_1^2 + 2a\omega_{3} \omega_2^2 + 2\omega_1u_2 ] - \omega_2\lambda_2(4a^2\omega_{3}^2\omega_1 - 3a\omega_{3}u_2 ) }{\dot{\lambda}_2 \omega_2 - \omega_{3}(1+2a)(\lambda_3x_2 - \lambda_4x_1)  }
\end{equation}
A similar result to that given in this derivation is obtained if $u_2$ is singular and $(u_1,u_3)$ bang-bang.  Due to this similarity, the derivation of the singular control for $u_2$ is not repeated here.  

\subsection{Derivation of Singular Control for \texorpdfstring{$\omega_3=0$}{TEXT}}\label{sect:appSpecial}

Consider the case where $u_1$ is singular and $u_2$ is bang-bang. The Hamiltonian for the reduced system of dynamic equations in Eq.~\eqref{eq:FS-eom} is
\begin{equation}\label{eq:reducedHam}
  \C{H}=\lambda_1u_1 + \lambda_2u_2 + \lambda_4 [\omega_2x_1x_2 + \frac{\omega_1}{2}(1+x_1^2-x_2^2))]+ \lambda_5 [\omega_1x_1x_2 + \frac{\omega_2}{2}(1+x_2^2-x_1^2)]
\end{equation}
On an interval where $u_1$ is singular, the switching function $g_1$ must be differentiated four times with respect to time until the control explicitly appears in order to determine the singular control.  First, it is known from Eq.~\eqref{eq:reducedHam} that
\begin{equation}
  g_1 = \lambda_1 = 0
\end{equation}
when $u_1$ is singular.  Differentiating $g_1$ with respect to time four times then yields the following results:
\begin{eqnarray}\label{eq:reduced-dg1dt}
  \frac{d g_1}{dt} & = &\lambda_{4}\frac{1+x_1^2-x_2^2}{2} - \lambda_{5}x_1x_2, \\\label{eq:reduced-d2g1dt2}
  \frac{d^2g_1}{dt^2} &  = &\omega_2(\lambda_4x_2-\lambda_5x_1), \\\label{eq:reduced-d3g1dt3}
  \frac{d^3g_1}{dt^3} & = &\omega_1\omega_2\dot{\lambda}_2, \\\label{eq:reduced-d4g1dt4}
  \frac{d^4g_1}{dt^4} & = & u_1\omega_2\dot{\lambda}_2,
\end{eqnarray}  
where the fact that $u_2$ is bang-bang and $u_2 \neq 0 \rightarrow \omega_2 \neq 0$ implies that $(\lambda_4x_2-\lambda_5x_1)=0$ in Eq.~\eqref{eq:reduced-d2g1dt2}. Given that $\omega_2 \neq 0$, if $\dot{\lambda}_2 \neq 0$ then Eq.~\eqref{eq:reduced-d3g1dt3} implies that $\omega_1=0$.
Note that all four derivatives in Eqs.~\eqref{eq:reduced-dg1dt}--\eqref{eq:reduced-d4g1dt4} are zero on a singular arc because $g_1$ is itself zero on the singular arc.  Therefore, $d^4g_1/dt^4$ in Eq.~\eqref{eq:reduced-d4g1dt4} is zero on a singular arc and the singular control $u_1$ is obtained as
\begin{equation}
  u_1 = 0.
\end{equation}
Now, the generalized Legendre-Clebsch condition implies that $\omega_2\dot{\lambda}_2 \geq 0$ for the singular control to be optimal which is equivalent to the condition $\omega_1=x_1=0$ along the singular arc (see Ref.~\cite{Shen1999} for details).

\renewcommand{\baselinestretch}{1.0}
\normalsize\normalfont
\bibliographystyle{aiaa}

\begin{thebibliography}{10}
\newcommand{\enquote}[1]{``#1''}

\bibitem{Shen1999}
Shen, H. and Tsiotras, P., \enquote{Time-Optimal Control of Axisymmetric Rigid
  Spacecraft Using Two Controls,} {\em Journal of Guidance, Control, and
  Dynamics\/}, Vol.~22, No.~5, Sept. 1999, pp.~682--694. \\
  \url{https://doi.org/10.2514/2.4436}.

\bibitem{Kirk2004}
Kirk, D.~E., {\em Optimal control theory: an introduction\/}, Courier
  Corporation, 2004.

\bibitem{Bryson1975}
Bryson, A.~E. and Ho, Y., {\em Applied Optimal Control: Optimization,
  Estimation, and Control\/}, Hemisphere Publishing Corporation, 1975.

\bibitem{Bilimoria1993}
Bilimoria, K.~D. and Wie, B., \enquote{Time-optimal three-axis reorientation of
  a rigid spacecraft,} {\em Journal of Guidance, Control, and Dynamics\/},
  Vol.~16, No.~3, 1993, pp.~446--452. \\ \url{https://doi.org/10.2514/3.21030}.

\bibitem{Tsiotras1995}
Tsiotras, P. and Longuski, J.~M., \enquote{A new parameterization of the
  attitude kinematics,} {\em Journal of the Astronautical Sciences\/}, Vol.~43,
  No.~3, 1995, pp.~243--262.

\bibitem{Fleming2008}
Fleming, A. and Ross, I.~M., \enquote{Optimal Control of Spinning Axisymmetric
  Spacecraft: A Pseudospectral Approach,} {\em {AIAA} Guidance, Navigation and
  Control Conference and Exhibit\/}, American Institute of Aeronautics and
  Astronautics, June 2008. \\ \url{https://doi.org/10.2514/6.2008-7164}.

\bibitem{Kim2013}
Kim, D., Turner, J.~D., and Leeghim, H., \enquote{Reorientation of Asymmetric
  Rigid Body Using Two Controls,} {\em Mathematical Problems in Engineering\/},
  Vol.~2013, 2013, pp.~1--8. \\ \url{https://doi.org/10.1155/2013/708935}.

\bibitem{Li2016}
Li, J., \enquote{Time-optimal three-axis reorientation of asymmetric rigid
  spacecraft via homotopic approach,} {\em Advances in Space Research\/},
  Vol.~57, No.~10, May 2016, pp.~2204--2217. \\
  \url{https://doi.org/10.1016/j.asr.2016.02.016}.

\bibitem{Lan2020}
Lan, J. and Li, J., \enquote{Newton-Kantorovich/Radau pseudospectral solution
  to rigid spacecraft time-optimal three-axis reorientation,} {\em Advances in
  Space Research\/}, Vol.~65, No.~11, June 2020, pp.~2662--2673. \\
  \url{https://doi.org/10.1016/j.asr.2020.02.023}.

\bibitem{Betts2010}
Betts, J.~T., {\em Practical methods for optimal control and estimation using
  nonlinear programming\/}, SIAM, 2010.

\bibitem{Gill2002}
Gill, P.~E., Murray, W., and Saunders, M.~A., \enquote{{SNOPT}: {A}n {SQP}
  {A}lgorithm for {L}arge-{S}cale {C}onstrained {O}ptimization,} {\em SIAM
  Review\/}, Vol.~47, No.~1, January 2002, pp.~99--131. \\
  \url{https://doi.org/10.1137/S0036144504446096}.

\bibitem{Biegler2008}
Biegler, L.~T. and Zavala, V.~M., \enquote{{L}arge-{S}cale {N}onlinear
  {P}rogramming {U}sing {IPOPT}: {A}n {I}ntegrating {F}ramework for
  {E}nterprise-{W}ide {O}ptimization,} {\em Computers and Chemical
  Engineering\/}, Vol.~33, No.~3, March 2008, pp.~575--582. \\
  \url{https://doi.org/10.1016/j.compchemeng.2008.08.006}.

\bibitem{Seywald1993}
Seywald, H. and Kumar, R.~R., \enquote{Singular Control in Minimum Time
  Spacecraft Reorientation,} {\em Journal of Guidance, Control, and
  Dynamics\/}, Vol.~16, No.~4, July 1993, pp.~686--694. \\
  \url{https://doi.org/10.2514/3.56607}.

\bibitem{Pager2022}
Pager, E.~R. and Rao, A.~V., \enquote{Method for solving bang-bang and singular
  optimal control problems using adaptive Radau collocation,} {\em
  Computational Optimization and Applications\/}, Jan. 2022. \\
  \url{https://doi.org/10.1007/s10589-022-00350-6}.

\bibitem{Kameswaran2008}
Kameswaran, S. and Biegler, L.~T., \enquote{{C}onvergence {R}ates for {D}irect
  {T}ranscription of {O}ptimal {C}ontrol {P}roblems {U}sing {C}ollocation at
  {R}adau {P}oints,} {\em Computational Optimization and Applications\/},
  Vol.~41, No.~1, 2008, pp.~81--126. \\
  \url{https://doi.org/10.1007/s10589--007--9098--9}.

\bibitem{Garg2010}
Garg, D., Patterson, M.~A., Hager, W.~W., Rao, A.~V., Benson, D.~A., and
  Huntington, G.~T., \enquote{{A} {U}nified {F}ramework for the {N}umerical
  {S}olution of {O}ptimal {C}ontrol {P}roblems {U}sing {P}seudospectral
  {M}ethods,} {\em Automatica\/}, Vol.~46, No.~11, November 2010,
  pp.~1843--1851. \\ \url{https://doi.org/10.1016/j.automatica.2010.06.048}.

\bibitem{Garg2011a}
Garg, D., Hager, W.~W., and Rao, A.~V., \enquote{{P}seudospectral {M}ethods for
  {S}olving {I}nfinite-{H}orizon {O}ptimal {C}ontrol {P}roblems,} {\em
  Automatica\/}, Vol.~47, No.~4, April 2011, pp.~829--837. \\
  \url{https://doi.org/10.1016/j.automatica.2011.01.085}.

\bibitem{Garg2011b}
Garg, D., Patterson, M.~A., Darby, C.~L., Francolin, C., Huntington, G.~T.,
  Hager, W.~W., and Rao, A.~V., \enquote{{D}irect {T}rajectory {O}ptimization
  and {C}ostate {E}stimation of {F}inite-{H}orizon and {I}nfinite-{H}orizon
  {O}ptimal {C}ontrol {P}roblems via a {R}adau {P}seudospectral {M}ethod,} {\em
  Computational Optimization and Applications\/}, Vol.~49, No.~2, June 2011,
  pp.~335--358. \\ \url{http://dx.doi.org/10.1007/s10589--009--9291--0}.

\bibitem{Patterson2015}
Patterson, M.~A., Hager, W.~W., and Rao, A.~V., \enquote{A $ph$ mesh refinement
  method for optimal control,} {\em Optimal Control Applications and
  Methods\/}, Vol.~36, No.~4, July--August 2015, pp.~398--421. \\
  \url{https://doi.org/10.1002/oca.2114}.

\bibitem{Hager2016}
Hager, W.~W., Hou, H., and Rao, A.~V., \enquote{Convergence Rate for a Gauss
  Collocation Method Applied to Unconstrained Optimal Control,} {\em Journal of
  Optimization Theory and Applications\/}, Vol.~169, No.~3, June 2016, pp.~801
  -- 824. \\ \url{https://doi.org/10.1007/s10957--016--0929--7}.

\bibitem{Hager2017}
Hager, W.~W., Hou, H., and Rao, A.~V., \enquote{Lebesgue Constants Arising in a
  Class of Collocation Methods,} {\em IMA Journal of Numerical Analysis\/},
  Vol.~37, No.~4, October 2017, pp.~1884--1901. \\
  \url{https://doi.org/10.1093/imanum/drw060}.

\bibitem{Hager2018}
Hager, W.~W., Liu, J., Mohapatra, S., Rao, A.~V., and Wang, X.-S.,
  \enquote{Convergence rate for a {Gauss} collocation method applied to
  constrained optimal control,} {\em SIAM Journal on Control and
  Optimization\/}, Vol.~56, 2018, pp.~1386--1411, \\
  \url{https://doi.org/10.1137/16M1096761}.

\bibitem{Hager2019}
Hager, W.~W., Hou, H., Mohapatra, S., Rao, A.~V., and Wang, X.-S.,
  \enquote{Convergence rate for a {Radau} {\it hp}-collocation method applied
  to constrained optimal control,} {\em Computational Optimization and
  Applications\/}, Vol.~74, 2019, pp.~274--314, \\
  \url{https://doi.org/10.1007/s10589--019--00100--1}.

\bibitem{Athans2013}
Athans, M. and Falb, P.~L., {\em Optimal control: an introduction to the theory
  and its applications\/}, Courier Corporation, 2013.

\bibitem{Schattler2012}
Sch{\"a}ttler, H. and Ledzewicz, U., {\em Geometric optimal control: theory,
  methods and examples\/}, Vol.~38, Springer Science \& Business Media, 2012.

\bibitem{Kelley1967}
Kelley, H., Kopp, R.~E., and Moyer, H.~G., \enquote{Topics in Optimization,
  edited by Leitman,} 1967.

\bibitem{Kopp1965}
Kopp, R.~E. and Moyer, H.~G., \enquote{Necessary conditions for singular
  extremals,} {\em {AIAA} Journal\/}, Vol.~3, No.~8, Aug. 1965, pp.~1439--1444.
  \\ \url{https://doi.org/10.2514/3.3165}.

\bibitem{Pager2022conf}
Pager, E.~R. and Rao, A.~V., \enquote{Structure Identification Method for
  Nonsmooth and Singular Optimal Control Problems,} {\em {AIAA} {SCITECH} 2022
  Forum\/}, American Institute of Aeronautics and Astronautics, Jan. 2022. \\
  \url{https://doi.org/10.2514/6.2022-1599}.

\bibitem{Weinstein2017}
Weinstein, M.~J. and Rao, A.~V., \enquote{Algorithm 984: ADiGator, a Toolbox
  for the Algorithmic Differentiation of Mathematical Functions in MATLAB Using
  Source Transformation via Operator Overloading,} {\em {ACM} Transactions on
  Mathematical Software\/}, Vol.~44, No.~2, Aug. 2017, pp.~1--25. \\
  \url{https://doi.org/10.1145/3104990}.

\bibitem{Patterson2014}
Patterson, M.~A. and Rao, A.~V., \enquote{{GPOPS}-{II},} {\em {ACM}
  Transactions on Mathematical Software\/}, Vol.~41, No.~1, oct 2014,
  pp.~1--37. \\ \url{https://doi.org/10.1145/2558904}.

\end{thebibliography}

\end{document}